\documentclass[11pt]{article}

\usepackage{amsmath,amsthm}

\usepackage{amssymb,latexsym}

\usepackage{enumerate}

\newtheorem{tw}{Theorem}[section]
\newtheorem{lm}[tw]{Lemma}
\newtheorem{wn}[tw]{Corollary}
\newtheorem{stw}[tw]{Proposition}
\newenvironment{dow}{\it Proof.\rm}{\hfill $\Box$}

\theoremstyle{definition}
\newtheorem*{df}{Definition}
\newtheorem{uw}[tw]{Remark}

\newcommand{\BF}{{\mathbb F}}

\newcommand{\BM}{{\mathbb M}}
\newcommand{\BN}{{\mathbb N}}
\newcommand{\BR}{{\mathbb R}}

\newcommand{\bM}{{\mathbf{M}}}
\newcommand{\bX}{{\mathbf{X}}}

\newcommand{\VV}{{\mathcal V}}
\newcommand{\WW}{{\mathcal W}}
\newcommand{\FF}{{\mathcal{F}}}

\newcommand{\BB}{{\mathcal{B}}}
\newcommand{\LL}{{\mathcal{L}}}

\newcommand{\MM}{{\mathcal{M}}}

\newcommand{\PP}{{\mathcal{P}}}
\newcommand{\EE}{{\mathcal{E}}}

\numberwithin{equation}{section}

\begin{document}
\title 
{Obstacle problem for
evolution equations involving measure data and operator
corresponding to semi-Dirichlet form}
\author{Tomasz Klimsiak}
\date{}
\maketitle
\begin{abstract}
In the paper, we consider the obstacle problem, with one and two irregular barriers, for semilinear evolution equation involving measure data and operator corresponding to a semi-Dirichlet form. We prove the existence and uniqueness of solutions under the assumption that the right-hand side of the equation is monotone and satisfies mild integrability conditions. To treat the case of irregular barriers, we extend the theory of precise versions of functions introduced by M. Pierre. We also give some applications to the so-called switching problem.
\end{abstract}


 \footnotetext{T. Klimsiak: Institute of Mathematics,
Polish Academy of Sciences, \'Sniadeckich 8, 00-956 Warszawa,
Poland, and Faculty of Mathematics and Computer Science, Nicolaus
Copernicus University, Chopina 12/18, 87-100 Toru\'n, Poland.
E-mail:
tomas@mat.umk.pl}

\section{Introduction}

Let $E$ be a locally compact separable metric space, $m$ be a
Radon measure on $E$ with full support, and let $\{B^{(t)},t\ge
0\}$ be a family of  regular semi-Dirichlet forms on $L^2(E; m)$
with common domain $V$ satisfying some  regularity conditions. By
$L_t$ denote the operator corresponding to the form $B^{(t)}$. In
the present paper we study the  obstacle problem with one and two
irregular barriers. In the case of one barrier $h:E\rightarrow\BR$
it can be stated as follows: for given $\varphi:E\rightarrow\BR$,
$f:[0,T]\times E\times\BR\rightarrow\BR$ and smooth (with respect to the
parabolic capacity Cap  associated with $\{B^{(t)},t\ge 0\}$)
measure $\mu$ on $E_{0,T}\equiv(0,T)\times E$ find
$u:\bar{E}_{0,T}\equiv (0,T]\times E\rightarrow\BR$ such that
\begin{equation}
\label{eq1.1} \left\{
\begin{array}{ll}-\frac{\partial u}{\partial t}-L_t u
= f(\cdot,u)+\mu&\mbox{on the set } \{u>h\},\quad
u(T,\cdot)=\varphi,
\smallskip\\
 -\frac{\partial u}{\partial t}-L_t u \ge f(\cdot,u)+\mu &\mbox{on  }
 E_{0,T},\smallskip\\
 u\ge h.
\end{array}
\right.
\end{equation}
In the second part of the paper we show how the results on
(\ref{eq1.1}) can be used to solve some system of variational
inequalities associated with  so-called switching problem.

Problems of the form  (\ref{eq1.1}) are  at present quite well
investigated in the case where $L_t$ are local operators.
Classical results for one or two regular barriers and
$L^2$-integrable data are to be found in the monograph \cite{BL}.
Semilinear obstacle problem with uniformly elliptic divergent form
operators $L_t$ and one or two irregular barriers was studied
carefully in \cite{K:SPA1,K:PA} in case of $L^2$-data, and in
\cite{KR:JEE} in case of measure data. In \cite{BDS} it is considered
evolutianry $p$-Laplacian type equation (with $p\in (1,\infty)$). In an important paper
\cite{Pierre1} linear problem (\ref{eq1.1}) with one irregular
barrier,  $L^2$-data  and operators $L_t$ associated with
Dirichlet  forms is  considered. The aim of the present paper is
to generalize or strengthen the existing results in the sense that  we
consider semilinear equations involving measure data  with two
barriers and  wide class of operators corresponding to
semi-Dirichlet forms. As for the obstacles, we only assume that
they are quasi-c\`adl\`ag  functions satisfying some mild integrability
conditions. The class of quasi-c\`adl\`ag functions naturally arises in the study of evolution equations. It includes quasi-continuous
functions and parabolic potentials (which in general are not
quasi-continuous). 

When considering problem (\ref{eq1.1}) with measure data, one of
the first difficulties one encounters is the proper definition of
a solution. Because of measure data, the usual variational
approach is not applicable. Moreover,  even in the case of
$L^2$-data, the variational inequalities approach does not give
uniqueness of solutions (see \cite{Mignot}). Therefore, following
\cite{Pierre1} and \cite{KR:JEE}, we consider  so-called
complementary system associated with (\ref{eq1.1}). Roughly
speaking, by a solution of this system we mean a pair $(u,\nu)$
consisting of a quasi-c\`adla\`g  function
$u:\bar{E}_{0,T}\rightarrow\BR$ and a positive smooth measure
$\nu$ on $E_{0,T}$ such that
\begin{equation}
\label{eq1.2.1} -\frac{\partial u}{\partial t}-L_t u =
f(\cdot,u)+\nu+\mu\quad\mbox{on } E_{0,T},\quad
u(T,\cdot)=\varphi\quad \mbox{in }E,
\end{equation}
\begin{equation}
\label{eq1.2.2}
 ``\nu\mbox{ is minimal}",
\end{equation}
\begin{equation}
\label{eq1.2.3} u\ge h\quad \mbox{q.e.},
\end{equation}
where q.e. means quasi-everywhere with respect to the capacity
Cap. Of course, in the above formulation one has to give rigorous
meaning to (\ref{eq1.2.1}) and (\ref{eq1.2.2}). As for
(\ref{eq1.2.1}), we develop some ideas from the paper \cite{K:JFA}
devoted to evolution equations involving measure data and
operators associated with semi-Dirichlet forms.

In the paper we assume that $\mu$ belongs to the class
$\mathbb{M}$ of smooth Borel measures with finite potential, which
 under additional assumption of duality for the family
$\{B^{(t)},t\ge 0\}$, takes the form
\[
\mathbb{M}=\bigcup_{\rho} \MM_\rho,
\]
where $\MM_\rho$ denotes the set of all smooth signed measures on
$E_{0,T}$ such that $\|\mu\|_\rho=\int_E\rho\,d|\mu|<\infty$, and
the union is taken over all strictly positive excessive functions
$\rho$.

If $\nu\in\mathbb{M}$, then we define a solution of
(\ref{eq1.2.1}) as in \cite{K:JFA}. To formulate this definition,
denote by $\bM$ a Hunt process $\{(\bX_t)_{t\ge 0},(P_z)_{z\in
E_{0,T}}\}$ with life time $\zeta$ associated with the operator
$\frac{\partial}{\partial t}+L_t$, and set $\zeta_\upsilon=
\zeta\wedge(T-\upsilon(0))$, where $\upsilon$ is the uniform
motion to the right, i.e. $\upsilon(t)=\upsilon(0)+t$ and
$\upsilon(0)=s$ under $P_z$ with $z=(s,x)$. By $A^{\mu},A^{\nu}$
denote natural additive functionals of $\bM$  in the Revuz
correspondence with $\mu$ and $\nu$, respectively. By a solution
of (\ref{eq1.2.1}) we mean a function $u$ such that for q.e. $z\in
E_{0,T}$,
\begin{equation}
\label{eq1.4} u(z)=E_z\Big(\varphi(\mathbf{X}_{\zeta_{\upsilon}})
+\int_0^{\zeta_\upsilon} f(\mathbf{X}_r, u(\mathbf{X}_r))\, dr
+\int_0^{\zeta_\upsilon}dA^\mu_r+\int_0^{\zeta_\upsilon}dA^\nu_r\Big).
\end{equation}
Formula (\ref{eq1.4}) may be viewed as a  nonlinear Feynman-Kac
formula.

Unfortunately, in general, the ``reaction measure" $\nu$ need not
belong to $\mathbb{M}$ (in fact, as shown in \cite{K:PA}, in case
of two barriers, it may be a nowhere Radon measure). In such case
we say that $u$ satisfies (\ref{eq1.2.1}) if $u$ is of class (D),
i.e. there is a potential on $\bar E_{0,T}$ (see Section 3.2) dominating
$|u|$ on $E_{0,T}$, and there exists a local martingale  $M$ (with $M_0=0$)
such that the following stochastic equation is satisfied
under the measure $P_z$ for q.e. $z\in E_{0,T}$:
\begin{align}
\label{eq1.3}
u(\mathbf{X}_t)&=\varphi(\mathbf{X}_{\zeta_{\upsilon}})
+\int_t^{\zeta_\upsilon} f(\mathbf{X}_r, u(\mathbf{X}_r))\, dr\nonumber\\
&\quad+\int_t^{\zeta_\upsilon}dA^\mu_r+\int_t^{\zeta_\upsilon}dA^\nu_r
-\int_t^{\zeta_\upsilon}\, dM_r,\quad t\in [0,\zeta_\upsilon].
\end{align}
In the above definition the requirement that $u$  is of class (D)
is very important. The reason is that (\ref{eq1.3}) is also
satisfied by functions solving equation (\ref{eq1.2.1}) with
additional nontrivial singular  (with respect to  Cap) measure on
its right-hand side  (see \cite{Kl:CVPDE}). These functions are
not of class (D).

If $\nu \in\mathbb{M}$ then (\ref{eq1.3}) is equivalent to
(\ref{eq1.4}). Furthermore, if the time dependent Dirichlet form
$\EE^{0,T}$ on $L^2(E_{0,T};m_1:=dt\otimes m)$ determined by the
family $\{B^{(t)},t\ge 0\}$ has the dual Markov property and
$\varphi\in L^1(E; m)$, $f(\cdot, u)\in L^1(E_{0,T};dt\otimes m)$,
$\mu,\nu\in\MM_1$ (i.e. $\mu,\nu$ are bounded), then the solution
$u$ in the sense of (\ref{eq1.4}) is a renormalized solution of
(\ref{eq1.2.1}) in the sense defined in \cite{KR:NDEA}. In
particular, this means that $u$ has some further regularity
properties and may be defined in purely analytical terms. More
precisely, $u$ is a renormalized solution if  the truncations
$T_ku=-k\vee (u\wedge k)$ of $u$ belong to the space $L^2(0,T;V)$,
and there exists a sequence $\{\lambda_{k}\}\subset\MM_{1}$ such
that $\|\lambda_{k}\|_{1}\rightarrow0$ as $k\rightarrow\infty$,
and for every bounded $v\in L^2(0,T;V)$ such that $\frac{\partial
v}{\partial t} \in L^2(0,T;V')$ and $v(0)=0$ we have
\begin{equation}
\label{eq1.5} \EE^{0,T}(T_{k}u,v) =(T_k\varphi,v(T))_{L^2} +(
f(\cdot, u), v)_{L^2}+\int_{E_{0,T}} v\, d\mu+\int_{E_{0,T}}v\,
d\nu_k
\end{equation}
for all $k\ge0$. In case of local operators of Leray-Lions type
the above definition of renormalized solutions was proposed in
\cite{PPP} (for the case of elliptic equations see \cite{DMOP}).

We now turn to condition (\ref{eq1.2.2}). Intuitively, it means
that $\nu$ ``acts only if necessary". If $h$ is quasi-continuous,
this statement means that $\nu$ acts only when $u=h$, because then
the right formulation of the minimality condition  takes the form
\begin{equation}
\label{eq1.6} \int_{E_{0,T}}(u-h)\,d\nu=0
\end{equation}
(see \cite{KR:JEE,Pierre1}). If $h$ is only quasi-c\`adl\`ag, the
situation is more subtle. For such $h$ the function $u$ satisfying
(\ref{eq1.4}) is also quasi-c\`adl\`ag, so the left-hand side of
(\ref{eq1.6}) is well defined, because $\nu$ is a smooth measure.
However, in general, the left-hand side is strictly positive. M.
Pierre has shown (in the case of linear equations with $L^2$-data
and Dirichlet forms), that for  general barrier $h$  the condition
\begin{equation}
\label{eq1.7} \int_{E_{0,T}} (\tilde{u}-\tilde{h})\,d\nu=0
\end{equation}
is satisfied. Here $\tilde{u}$ is a precise $m_1$-version of $u$
and $\tilde{h}$ is an associated precise version of $h$. The
notions of a precise $m_1$-version of a potential and an
associated precise version of a function $h$ dominated by a
potential (which is not necessarily $m_1$-version of $h$) were
introduced in \cite{Pierre1,Pierre2}. In the paper, in case of
semi-Dirichlet forms, we use probabilistic methods to introduce
another notion of a precise $m_1$-version $\hat{u}$ of a
quasi-c\`adl\`ag function $u$ (note that potentials are
quasi-c\`adl\`ag). Since our barriers as well as solutions to
(\ref{eq1.2.1})--(\ref{eq1.2.3}) are quasi-c\`adl\`ag, we do not
need the notion of an associated precise version. We show that if
$u$ is an $L^2$ potential, then
\[
\hat{u}=\tilde{u}\quad \mbox{q.e},
\]
and for any quasi-c\`adl\`ag function $u$  which is dominated by
an $L^2$ potential,
\[
\hat{u}\le \tilde{u}\quad \mbox{q.e.},\quad
\int(\tilde{u}-\hat{u})\,d\nu=0.
\]
It follows in particular that in case of $L^2$ data and Dirichlet
form, (\ref{eq1.7}) is equivalent to
\begin{equation}
\label{eq1.7i} \int_{E_{0,T}} (\hat{u}-\hat{h})\,d\nu=0
\end{equation}

One reason why we introduce a new notion of a precise version is
that it is applicable to wider classes of operators and functions
then those considered in  \cite{Pierre1,Pierre2}. The second is
that our definition is more direct that the construction of a
precise version given in \cite{Pierre1,Pierre2}.
Namely, in our approach by a precise version of a quasi-c\`adl\`ag
function $u$ on $\bar{E}_{0,T}$  we mean a function $\hat u$ on
$\bar{E}_{0,T}$ such that for q.e. $z\in E_{0,T}$,
\[
\hat{u}(\bX_{t-})=u(\bX)_{t-}\,,\quad t\in (0,\zeta_\upsilon)
\]
As a consequence,
our definition appears to be very convenient for studying obstacle
problems and  may be applied to quite general  class of equations
(possibly nonlinear with measure data and two obstacles).

In case of one obstacle, the main result of the paper says that if
$\varphi,f(\cdot,0)$ satisfy mild integrability conditions, and
$f$ is monotone and continuous with respect to $u$, then for every
$\mu\in\mathbb{M}$ there exists a unique solution $(u,\nu)$ of
(\ref{eq1.1}), i.e.  a unique pair $(u,\nu)$ consisting of a
quasi-c\`adl\`ag function $u$ on $\bar{E}_{0,T}$ and positive
smooth measure $\nu$ on $E_{0,T}$ such that (\ref{eq1.3}),
(\ref{eq1.7i}) and (\ref{eq1.2.3}) are satisfied. We also give
conditions under which $\nu\in\mathbb{M}$ or $\nu\in \MM_1$, i.e.
when equivalent to (\ref{eq1.3}) formulations (\ref{eq1.4}) or
(\ref{eq1.5}) may be used. Moreover, we show that $u_n\nearrow u$
q.e., where $u_n$ is a solution of the following Cauchy problem
\begin{equation}
\label{eq1.8} -\frac{\partial u_n}{\partial t}-L_t u_n =
f(\cdot,u_n)+n(u_n-h)^-+\mu,\quad u_n(T,\cdot)=\varphi.
\end{equation}

Our probabilistic approach allows us to prove  similar results
also for two quasi-c\`adl\`ag barriers $h_1,h_2$ satisfying some
separation condition. In the case of two barriers, the measure
$\nu$ appearing in the definition of a solution is a signed smooth
measure. Moreover, we replace the minimality  condition
(\ref{eq1.7}) by
\begin{equation}
\label{eq1.9} \int_{E_{0,T}}
(\hat{u}-\hat{h}_1)\,d\nu^+=\int_{E_{0,T}}(\hat{h}_2-\hat{u})\,d\nu^-=0,
\end{equation}
and replace condition (\ref{eq1.2.3}) by  $h_1\le u\le h_2$ q.e.
We show that under the same conditions on $\varphi,f,\mu$  as in
the case of one barrier, there exists a unique solution $(u,\nu)$
of the obstacle problem with two barriers. We also  show that
$\bar{u}_n\nearrow u$ and $\lambda_n\nearrow \nu^-$, where
$(\bar{u}_n,\lambda_n)$ is a solution of  problem of the form
(\ref{eq1.1}), but with one upper barrier $h_2$ and $f$ replaced
by
\[
f_n(t,x,y)=f(t,x,y)+n(y-h_1(t,x))^{-}.
\]
We prove these results under two different separation conditions.
The first one, more general,  may be viewed as some analytical
version of the  Mokobodzki condition considered in the literature
devoted to reflected stochastic differential equations (see, e.g.,
\cite{K:SPA2}). The second one, which is simpler and usually easier
to check, as yet, has not been considered in the literature on
evolution equations. It says that
\begin{equation}
\label{eq1.13} h_1<h_2,\quad \hat{h}_1<\hat{h}_2\quad\mbox{q.e.}
\end{equation}
If $h_1, h_2$ are quasi-continuous, then (\ref{eq1.13}) reduces to
the condition $h_1<h_2$ q.e., because  quasi-continuous functions
are precise (see \cite{K:PA} for this case).

Note also that at the end of the Section 6 we show that the study of the obstacle problem with one merely measurable barrier (or two measurable barriers satisfying the Mokobodzki condition) can be reduced to the study of the obstacle problem with quasi-c\`adl\`ag barriers. It should be stressed, however, that when dealing with merely measurable barriers, our definition of a solution is weaker. Namely, instead of (\ref{eq1.2.3}) we only require that $u\ge h$ $m_1$-a.e. (or $h_1\le u\le h_2$ $m_1$-a.e. in case of two barriers).

In the last part of the paper  we use our results on (\ref{eq1.1})
to study  so-called switching problem (see Section \ref{sec7}).
This problem is closely related to  system of quasi-variational
inequalities, which when written as a complementary system, has
the form
\begin{equation}
\label{eq1.13.1} -\frac{\partial u^j}{\partial t}-L_tu^{j}
=f^{j}(t,x,u)+\nu^{j}+\mu^j,
\end{equation}
\begin{equation}
\label{eq1.13.2} \int_{E_{0,T}}(\hat{u}^{j}-\hat
H^j(\cdot,u))\,d\nu^{j}=0,
\end{equation}
\begin{equation}
\label{eq1.13.3} u^{j}\ge H^j(\cdot,u)\quad\mbox{q.e.},
\end{equation}
where
\[
H^j(z,y)=\max_{i\in A_j} h_{j,i}(z,y^i),\quad z\in
E_{0,T},\,y\in\BR^N.
\]
In (\ref{eq1.13.1})--(\ref{eq1.13.3}), we are given
$f^j:E_{0,T}\times \BR^N\rightarrow\BR^N,\, h_{j,i}:
E_{0,T}\times\BR\rightarrow\BR$, $\mu^j\in\mathbb{M}$ and sets
$A_j\subset \{1,\dots, j-1,j+1,\dots,N\}$, and we are looking for
a pair $(u,\nu)= ((u^1,\dots,u^N),(\nu^1,\dots,\nu^N))$ satisfying
(\ref{eq1.13.1})--(\ref{eq1.13.3}) for $j=1,\dots,N$. Note that in
(\ref{eq1.13.3}) the barrier $H^j$ depends on $u$.

Systems of the form (\ref{eq1.13.1})--(\ref{eq1.13.3}) were
subject to numerous investigations, but only in the framework of
viscosity solutions and for  local operators  (see
\cite{DHP,EH,HJ,HM,HT,HZ,LNO1}) or for special class of nonlocal
operators associated with a random Poisson measure (see \cite{LNO}).
In the paper we prove an  existence result for
(\ref{eq1.13.1})--(\ref{eq1.13.3}). In the important special case
where
\[
h_{j,i}(z,y)= -c_{j,i}(z)+y^i,
\]
we show that there is a unique solution of
(\ref{eq1.13.1})--(\ref{eq1.13.3}), and moreover, that $u$  is
the value function for the optimal switching problem related to
(\ref{eq1.13.1})--(\ref{eq1.13.3}).

\section{Preliminaries}
\label{sec2}

In the paper $E$ is a locally compact separable metric space and
$m$ a Radon measure on $E$ such that supp$[m]=E$. For $T>0$ we set
$E_{0,T}=(0,T)\times E$, $\bar E_{0,T}=(0,T]\times E$.

Recall (see \cite{Oshima}) that a form $(B,V)$ is called
semi-Dirichlet on $L^2(E;m)$ if $V$ is a dense linear subspace of
$L^2(E;m)$, $B$ is a bilinear form on $V\times V$, and the
following conditions (B1)--(B4) are satisfied:
\begin{enumerate}
\item[(B1)]  $B$ is lower bounded, i.e. there exists $\alpha_0\ge 0$
such that $B_{\alpha_0}(u,u)\ge 0$ for $u\in V$, where
$B_{\alpha_0}(u,v)=B(u,v)+\alpha_0 (u,v)$,
\item[(B2)] $B$ satisfies the sector condition, i.e.  there exists $K>0$ such that
\[
|B(u,v)|\le KB_{\alpha_0}(u,u)^{1/2}B_{\alpha_0}(v,v)^{1/2},\quad
u,v\in V,
\]
\item[(B3)] $B$ is closed, i.e. for every $\alpha>\alpha_0$
the space $V$ equipped with the inner product
$B^{(s)}_{\alpha}(u,v)=\frac12(B_\alpha(u,v)+B_\alpha(v,u))$ is a
Hilbert space,
\item[(B4)] $B$ has the Markov property, i.e. for every $a\ge
0$, $ B(u\wedge a,u\wedge a)\le B(u\wedge a,u)$ for all $u\in V$.
\end{enumerate}

Condition (B4) is called the Markov property, because it is
equivalent to the fact that the semigroup $\{T_t,t\ge 0\}$
associated with $(B,V)$  is sub-Markov (see \cite[Theorem
1.1.5]{Oshima}). Recall that $(B,V)$ is said to have the dual
Markov property if
\begin{enumerate}
\item[(B5)] for every $a\ge 0$,
\[
B(u\wedge a,u\wedge a)\le B(u,u\wedge a),\quad u\in V.
\]
\end{enumerate}
Note that (B5)  is equivalent to the fact that the dual semigroup
$\{\hat{T}_t,t\ge 0\}$ associated with $(B,V)$ is sub-Markov (see
\cite[Theorem 1.1.5]{Oshima}). For the notions of transiency and
regularity see \cite[Sections 1.2, 1.3]{Oshima}.

In the paper  $\{B^{(t)},t\in\BR\}$ is  a family of regular
semi-Dirichlet forms  on $L^2(E;m)$ satisfying (B2), (B3) with
some constants $K,\alpha_0$ not depending on $t$. We also assume
that for all $u,v\in V$ the mapping $\BR\ni t\mapsto B^{(t)}(u,v)$
is measurable, and for some $\lambda\ge1$,
\[
\lambda^{-1} B^{(0)}(u,u)\le B^{(t)}(u,u)\le \lambda
B^{(0)}(u,u),\quad u\in V,\quad t\in\BR.
\]

Let $(\EE, D[\EE])$ denote the time-dependent semi-Dirichlet form on $L^2(E_1;m_1)$
($E_1:=\BR\times E$)
associated with the family $\{B^{(t)},t\in\BR\}$ (see
\cite[Section 6.1]{Oshima}), and Cap denote the associated
capacity. We say that some property holds quasi-everywhere (q.e.
for short) if it holds outside some set $B\subset E_1$ such that Cap$(B)=0$.
Capacity Cap on $E_{0,T}$ is equivalent to the capacity considered in \cite{Pierre1,Pierre2}
(see \cite{Pierre3}) in the context of parabolic variational inequalities.

Let $\mu$ be a signed measure on $E_1$. By $|\mu|$ we
denote the variation of $\mu$, i.e. $|\mu|=\mu^{+}+\mu^{-}$, where
$\mu^{+}$ (resp. $\mu^{-}$) denote the positive (resp. negative)
part of $\mu$. Recall that a Borel measure on $E_1$ is
called smooth if $\mu$ charges no exceptional sets and there
exists an increasing sequence $\{F_n\}$ of closed subsets of
$E_1$ such that $|\mu|(F_n)<\infty$ for $n\ge1$ and
${\rm{Cap}}(K\setminus F_n )\rightarrow 0 $ for every compact
$K\subset E_1$.

It is known (see \cite[Section 6.3]{Oshima}) that there exists a
unique Hunt process $\bM=\{(\bX)_{t\ge0},(P_z)_{z\in
E_1}\}$ with life time $\zeta$ associated with the form $(\EE,
D[\EE])$. Moreover,
\[
\mathbf{X}_t=(\upsilon(t),X_{\upsilon(t)}),
\]
where $\upsilon$ is the uniform motion to the right, i.e.
$\upsilon(t)=\upsilon(0)+t$ and $\upsilon(0)=s,\, P_z$-a.s. for
$z=(s,x)$ (see \cite[Theorem 6.3.1]{Oshima}). In the sequel, for
fixed $T>0$ we set
\[
\zeta_\upsilon=\zeta\wedge(T-\upsilon(0)).
\]
Let us recall that from \cite[Lemma 6.3.2]{Oshima} it follows that a nearly Borel set $B\subset E_{1}$
is of capacity zero iff it is exceptional, i.e.  $P_z(\exists _{t>0};\,\, \bX_t\in B)=0$, q.e.
Recall also (see \cite[Section 6]{Oshima}) that there is one-to-one
correspondence, called Revuz duality, between positive smooth
measures and  positive natural additive functionals (PNAFs) of
$\bM$. For a positive smooth measure $\mu$ we denote by $A^\mu$
the unique PNAF in the Revuz duality with $\mu$. For a signed
smooth measure $\mu$ we put $A^\mu= A^{\mu^+}-A^{\mu^-}$. For a
fixed positive measurable function $f$ and a positive Borel
measure $\mu$ we denote by $f\cdot \mu$ the measure defined as
\[
(f\cdot\mu)(\eta)=\int_{E_1} \eta f\,d\mu,\quad
\eta\in\mathcal{B}^+(E_1).
\]

By $S_{0}$ we denote the set of all measures of finite energy
integrals, i.e. the set of all smooth measures $\mu$ having the
property that there is $K\ge0$ such that
\[
\int_{E}|\tilde{\eta}|\,d|\mu|\le K\|\eta\|_{\WW}, \quad \eta\in\WW,
\]
where $\tilde{\eta}$ is a quasi-continuous $m_1$-version of $\eta$
(for the existence of such version see \cite[Theorem
6.2.11]{Oshima}). By $\mathbb M$ we denote the set of all smooth
measures on $E_{0,T}$ such that
$E_zA^{|\mu|}_{\zeta_\upsilon}<\infty$ for q.e. $z\in E_{0,T}$.
$\mathbb{M}_c$  is the set of those $\mu\in\BM$  for which  $A^\mu$ is continuous.

For a given positive smooth measure $\mu$ on $\bar E_{0,T}$ we set
\[
R^{0,T}\mu(z)=E_zA^\mu_{\zeta_\upsilon},\quad z\in E_{0,T}.
\]
Set
\[
\VV_{0,T}=L^2(0,T; V),\qquad \WW_{0,T}=\{u\in\VV_{0,T}:
\frac{\partial u}{\partial t}\in \VV_{0,T}'\},
\]
\[
\WW_T=\{u\in\WW_{0,T}:u(T)=0\},\qquad
\WW_0=\{u\in\WW_{0,T}:u(0)=0\}
\]
and
\[
\EE^{0,T}(u,v)=\left\{
\begin{array}{ll}\int^T_0
\langle-\frac{\partial u}{\partial t}(t),v(t)\rangle\,dt
+\int_0^TB^{(t)}(u(t),v(t))\,dt,\quad
(u,v)\in\WW_T\times\VV_{0,T},
\medskip \\
\int^T_0\langle u(t),\frac{\partial v}{\partial t}(t)\rangle\,dt +
\int_0^TB^{(t)}(u(t),v(t))\,dt,\quad (u,v)\in\VV_{0,T}\times\WW_0,
\end{array}
\right.
\]
where $\langle\cdot,\cdot\rangle$ is  the duality pairing between
$V'$ and $V$ ($V'$ stands for the dual of $V$). It is known (see
\cite[Example I.4.9(iii)]{Stannat}) that $\EE^{0,T}$ is a
generalized semi-Dirichlet form. The operator associated with
$\EE^{0,T}$ has the form
\[
\LL=-\frac{\partial}{\partial t}-L_{t}, \quad D(\LL)=\{u\in
\WW_T:\LL u\in L^2(E_{0,T}; m_1) \},
\]
where $(L_t,D(D_t))$ is the operator associated with
$(B^{(t)},V)$, $t\in[0,T]$. By $(G^{0,T}_{\alpha})_{\alpha>0}$ we
denote the (unique) strongly continuous contraction resolvent on
$L^2(0,T;L^2(E;m))$  associated with $\EE^{0,T}$ (see Propositions
I.3.4 and I.3.6 in \cite{Stannat}).

\section{Precise versions of quasi-c\`adl\`ag functions}
\label{sec4}

\subsection{Precise versions  of parabolic potentials in the sense of Pierre and its probabilistic interpretation}

We first recall the notion of a precise version of a parabolic
potential introduced in \cite{Pierre2}.

\begin{df}
A measurable function $u\in\VV_{0,T}\cap L^\infty(0, T; L^2(E;m))$
is called a parabolic potential if for every nonnegative
$v\in\WW_0$,
\[
\int_0^T\Big\langle\frac{\partial v}{\partial t}(t),
u(t)\Big\rangle\,dt+ \int_0^TB^{(t)}(u(t), v(t))\,dt\ge 0.
\]
\end{df}
The set of all parabolic potentials will be denoted by $\PP^2$.

\begin{stw}
\label{stw3.1} Let $u\in\PP^2$. Then there exists a unique positive $\mu\in
S_0$ on $\bar E_{0,T}$ such that
\begin{equation}
\label{eq3.1} u(z)=E_zA_{\zeta_\upsilon}^\mu
\end{equation}
for $m_1$-a.e. $z\in E_{0,T}$.
\end{stw}
\begin{dow}
By \cite[Theorem  III.1]{Pierre2} there exists a positive measure
$\mu\in S_0$ on $\bar E_{0,T}$ such that
\[
\int_t^T(\frac{\partial v}{\partial t}(s), u(s))\,ds+
\int_t^TB^{(s)}(u(s),v(s))\,ds=\int_{(t,T]\times E}v\,d\mu- (v(t),
u(t))_{L^2}
\]
for all $v\in\WW_{0,T}$ and $t\in[0, T]$. This and \cite[Theorem
3.7]{K:JFA} yield (\ref{eq3.1}).
\end{dow}
\medskip

In the sequel, for $u\in\PP^2$ we set $\EE(u)=\mu$, where $\mu\in
S_0$ is the  measure from Proposition \ref{stw3.1}.

\begin{df}
A measurable function $u$ on $E_{0,T}$ is called precise (in the
sense of M. Pierre) if there exists a sequence $\{u_n\}$ of
quasi-continuous parabolic potentials such that $u_n\searrow u$
q.e. on $E_{0,T}$.
\end{df}

The following result has been proved in \cite{Pierre2}.

\begin{tw}
Each $u\in\PP^2$ has a precise $m_1$-version.
\end{tw}

In what follows we denote by $\tilde{u}$ a precise version of
$u\in\PP^2$ in the sense of Pierre. It is clear that $\tilde u$ it is defined q.e. In
\cite{Pierre2} it is proved that the mapping $t\mapsto\tilde
u(t,\cdot)\in L^2(E;m)$ is left continuous and has right limits,
whereas in  \cite[Proposition 3.4]{K:JFA} it is proved that $u$ defined by
(\ref{eq3.1}) has the property that the mapping $t\mapsto
u(t,\cdot)\in L^2(E;m)$ is right continuous and has left limits.
Since for any $B\in\BB(E)$ and $t\in[0, T]$, Cap$(\{t\}\times
B)=0$ if and only if $m(B)=0$ (see \cite[Proposition II.4]{Pierre2}), it follows that in general
Cap$(\{u\neq\tilde{u}\})>0$. In the sequel, for given  $u\in\PP^2$
we will always consider its version defined by (\ref{eq3.1}).

Let us recall that function $x:[a,b]\rightarrow \BR$ is called c\`adl\`ag (resp. c\`agl\`ad)
iff $x$ is right-continuous (resp. left-continuous) and has left (resp. right) limits.

\begin{lm}
\label{lm3.1} Assume that $\{x_n\}$ is a decreasing sequence of
c\`agl\`ad functions on $[0,T]$ such that $x_n(t)\searrow x(t)$,
$t\in[0, T]$, and $x(t)=-a(t)+b(t)$, $t\in[0, T]$,  for some
nondecreasing function $a$ and c\`agl\`ad function $b$ on $[0,T]$.
Then $a$ and $x$ are c\`agl\`ad functions.
\end{lm}
\begin{dow}
The proof is analogous to that of \cite[Lemma 2.2]{Peng}, so we
omit it.
\end{dow}

\begin{lm}
\label{lm3.2} Assume that for each $n\in\BN$,
\begin{equation}
Y_t^n=Y_0^n-A_t^n+M_t^n, \quad t\in[0, T],
\end{equation}
where $A^n$ is a predictable increasing  process
with $A^n_0=0$, and $M^n$ is a local  martingale with
$M^n_0=0$. If $Y^n$ is positive, $Y_t^n\searrow Y_t$, $t\in[0, T]$, and $E\sup_{t\le
T }(|Y^1_t|^2+|Y_t|^2)<\infty$, then
$\hat{Y}_t:=\lim_{n\rightarrow\infty} Y^n_{t-}$, $t\in[0, T]$, is
a c\`agl\`ad process.
\end{lm}

\begin{dow}
By \cite[Theorem 3.1]{LLP},
\[
E|A^n|_T^2+E[M^n]_T\le c E\sup_{t\in[0, T]}|Y^n_t|^2\le
cE\sup_{t\in[0, T]}(|Y^1_t|^2\vee|Y_t|^2).
\]
In particular, $\sup_nE|M^n_T|^2<\infty$. Therefore there exists
$X\in L^2(\FF_T)$ such that $M^n_T\rightarrow X$ weakly in
$L^2(\FF_T)$. Let $N$ be a c\`adl\`ag version of $E(X|\FF_t)$,
$t\in[0, T]$. Then for every $\tau\in\mathcal{T}$,
$M^n_\tau\rightarrow N_\tau$ weakly in $L^2(\FF_T)$. Indeed, if
$Z\in L^2(\FF_T)$, then
\begin{align}
\label{eq3.8} EM^n_\tau Z=E(E(M^n_T|\FF_\tau)\cdot Z)
&=E(M^n_TE(Z|\FF_\tau))\nonumber\\
&\rightarrow
EXE(Z|\FF_\tau)=E(E(X|\FF_\tau)Z)=EN_\tau\cdot Z.
\end{align}
Since $\bM$ is a Hunt process each
$\FF$-martingale $M$ has the property that $M_{\tau-}=M_\tau$ for
all predictable $\tau\in\mathcal{T}$ (see \cite[Proposition 2]{CW}). Therefore  for any
predictable $\tau\in\mathcal{T}$,
\begin{equation}
\label{eq3.9} A^n_{\tau-}=Y_0^n-Y^n_{\tau-}+M_{\tau-}^n=Y_0^n-
Y^n_{\tau-}+M_\tau^n.
\end{equation}
Set
\[
\hat{A}_t=Y_0-\hat{Y}_t+N_{t-}, \quad t\in[0, T].
\]
Of course, $\hat{A}$ is predictable. By (\ref{eq3.8}) and
(\ref{eq3.9}),  for any predictable $\tau\in\mathcal{T}$,
\[
A^n_{\tau-}\rightarrow Y_0-\hat{Y}_\tau+N_\tau=
Y_0-\hat{Y}_\tau+N_{\tau-}=\hat{A}_\tau
\]
weakly in $L^2(\FF_T)$. From the above convergence,
$\hat{A}_\sigma\le\hat{A}_\tau$ for all predictable $\sigma,
\tau\in\mathcal{T}$ such that $\sigma\le\tau$. Therefore applying
the predictable cross-section theorem (see \cite[[Theorem 86, p.
138]{DM}) we conclude that $\hat{A}$ is an increasing process.
Consequently, by Lemma \ref{lm3.1}, $\hat{Y}$ is c\`agl\`ad.
\end{dow}

Now we are ready to prove the main result of this section. In the
sequel, for a function $u$ on $\bar E_{0,T}$ we set
$u(\bX)_{t-}=Y_{t-}$ and $u(\bX)_{t+}=Y_{t+}$, where
$Y_t=u(\bX_t)$, $t\in(0,T]$.

\begin{tw}
\label{tw3.1} Assume that $u\in\PP^2$. Then for q.e. $z\in
E_{0,T}$,
\begin{equation}
\label{eq3.6} \tilde{u}(\bX_{t-})=u(\bX)_{t-}\,, \quad t\in(0,
\zeta_\upsilon], \quad P_z\mbox{\rm{-a.s.}}
\end{equation}
\end{tw}
\begin{dow}
By the definition of a precise version of a potential, there
exists a sequence $\{u_n\}$ of quasi-continuous parabolic
potentials such that $u_n\searrow \tilde{u}$ q.e. on  $E_{0, T}$.
Since $\varphi(u)\in \PP^2$ for  any bounded concave function $\varphi$ on $\BR^+$ and any $u\in \PP^2$ we may assume
that $u_n, u$ are bounded. Since $u_n, u\in\PP^2$,  then by Proposition \ref{stw3.1} and the
strong Markov property there exist measures $\mu_n, \mu\in\mathbb{M}$ and  martingales $M^n, M$ such that
\[
u_n(\bX_t)=\int_t^{\zeta_\tau}dA_r^{\mu_n}
-\int_t^{\zeta_\tau}dM_r^n,
\quad t\in[0, \zeta_\tau]
\]
and
\begin{equation}
\label{eq3.3}
u(\bX_t)=\int_t^{\zeta_\tau}dA_r^{\mu}-\int_t^{\zeta_\tau}dM_r,
\quad t\in[0, \zeta_\tau].
\end{equation}
Since $u_n\searrow\tilde{u}$
q.e. on $E_{0, T}$, we have
\[
u_n(\bX_{t-})\searrow\tilde{u}(\bX_{t-}), \quad t\in(0, \zeta_\tau],
\quad P_z\mbox{-a.s.}
\]
for q.e. $z\in E_{0, T}$. But $u_n(\bX_{t-})=u_n(\bX)_{t-}$,
$t\in(0, \zeta_{\tau}]$, $P_z$-a.s. for q.e. $z\in E_{0, T}$,
because $u_n$ is quasi-continuous. Therefore applying Lemma
\ref{lm3.2} shows that $\tilde{u}(\bX_{t-})$ is c\`agl\`ad. We are
going to show that
\begin{equation}
\label{eq3.4} \tilde{u}(\bX_-)_{t+}=u(\bX_t), \quad t\in(0,
\zeta_\tau).
\end{equation}
Since $\tilde{u}=u$ $m_1$-a.e.,
\begin{align*}
0=R^{0,T}|\tilde{u}-u|(z)&=E_z\int_0^{\zeta_\upsilon}|\tilde{u}(\bX_t)-u(\bX_t)|\,dt\\&=
E_z\int_0^{\zeta_\upsilon}|\tilde{u}(\bX_{t-})-u(\bX_t)|\,dt=
E_z\int_0^{\zeta_\upsilon}|(\tilde{u}(\bX_-)_{t+}-u(\bX_t)|\,dt
\end{align*}
for q.e. $z\in E_{0, T}$. Hence
\begin{equation}
\label{eq3.06} \tilde{u}(\bX_-)_{t+}=u(\bX_t)\quad\mbox{for a.e.
}t\in (0, \zeta_\tau)
\end{equation}
for q.e. $z\in E_{0, T}$. Since $\tilde{u}(\bX_-)_{t+}$ and
$u(\bX_t)$ are c\`adl\`ag processes, (\ref{eq3.06}) implies
(\ref{eq3.4}). In turn, (\ref{eq3.4}) implies (\ref{eq3.6}),
because $\tilde{u}(\bX_-)$ is c\`agl\`ad.
\end{dow}

\subsection{Probabilistic approach to precise versions
of quasi-c\`adl\`ag functions}

\begin{df}
We say that $u:\bar{E}_{0,T}\rightarrow\BR$ is quasi-c\`adl\`ag if
the process $u(\mathbf{X})$ is c\`adl\`ag on $[0,\zeta_\upsilon]$
under the measure $P_z$ for q.e. $z\in E_{0,T}$.
\end{df}

Of course, any quasi-continuous function is quasi-c\`adl\`ag.
In the sequel the class of smooth measures $\mu$
on $\bar E_{0,T}$ for which $E_zA^{|\mu|}_{\zeta_\upsilon}<\infty$ q.e. on $E_{0,T}$ we denote by $\mathbb M_T$. We say that $u$ is a  potential on $\bar E_{0,T}$ iff  for a positive $\mu\in \mathbb M_T$.
\begin{equation*}
u(z)=E_zA^\mu_{\zeta_\upsilon}
\end{equation*}
for q.e. $z\in E_{0,T}$.  By \cite[Proposition
3.4]{K:JFA}, each  potential on $\bar E_{0,T}$ is quasi-c\`adl\`ag.
By Proposition \ref{stw3.1}, each $u\in\PP^2$ is a  potential on $\bar E_{0,T}$.

\begin{tw}
\label{tw4.1} Let $u$ be a quasi-c\`adl\`ag function on
$\bar{E}_{0,T}$. Then there exists a unique (q.e.) function $\hat{u}$ such that
for q.e. $z\in E_{0,T}$,
\begin{equation}
\label{eq4.1} \hat{u}(\bX_{t-})=u(\bX)_{t-}\,,\quad t\in
(0,\zeta_\upsilon),\quad
P_z\mbox{-a.s.}
\end{equation}
\end{tw}
\begin{dow}
By \cite[Theorem 3.3.6]{Oshima} there exists a measure $\hat{m}_1$
equivalent to $m_1$ such that $(\mathbf{X},P_z)$ has the dual Hunt
process $(\hat{\mathbf{X}},\hat{P}_z)$ with respect to
$\hat{m}_1$. Therefore by \cite[Theorem 16.4]{GS} applied to the
process $u(X)_{-}$ there exists a Borel measurable function
$\hat{u}$ satisfying (\ref{eq4.1}). Uniqueness follows from the
very definition of exceptional sets and \cite[Proposition
11.1]{GS}.
\end{dow}
\medskip

\begin{df}
For a quasi-c\`adl\`ag function $u$ on $\bar E_{0,T}$ we call the
function $\hat{u}$ from  Theorem \ref{tw4.1} a precise
$m_1$-version of $u$.
\end{df}

\begin{wn}
For $u\in\PP^2$, $\tilde{u}=\hat{u}$ q.e.
\end{wn}

\begin{uw}
\label{uw3.7} It is clear that $\hat{u}$ is  an $m_1$-version of
$u$ (see the reasoning before (\ref{eq3.06})).
\end{uw}

In the sequel we will need the following result.

\begin{stw}
\label{stw4.1} Let $\mu$ be a positive smooth measure on
$E_{0,T}$, and let $u$  be a positive quasi-c\`adl\`ag function on
$E_{0,T}$. Then
\begin{equation}
\label{eq3.7}
\int_0^{\zeta_\upsilon}[u(\mathbf{X})]_{t-}\,dA^{\mu}_t
=\int_0^{\zeta_\upsilon}\hat{u}(\mathbf{X}_t)\,dA^{\mu}_t,\quad
P_z\mbox{-a.s.}
\end{equation}
for q.e. $z\in E_{0,T}$. Moreover,
\[
E_z\int_0^{\zeta_\upsilon}u(\mathbf{X}_t)\,dA^{\mu}_t =0\quad
\mbox{for q.e.}\quad z\in E_{0,T}
\]
if and only if $\int_{E_{0,T}}u\,d\mu=0$.
\end{stw}
\begin{dow}
By Theorem \ref{tw4.1},
$\int_0^{\zeta_\upsilon}[u(\mathbf{X})]_{t-}\,dA^{\mu}_t
=\int_0^{\zeta_\upsilon}\hat{u}(\mathbf{X}_{t-})\,dA^{\mu}_t$,
$P_z$-a.s. It is well known that $\mathbf{X}$ has only totally
inaccessible  jumps (as a Hunt process), which combined with the fact that $A^\mu$ is
predictable gives (\ref{eq3.7}). The second part of the
proposition follows directly from the Revuz duality.
\end{dow}

\subsection{Associated precise versions in the sense of Pierre}

Let $u$ be a measurable function on $\bar E_{0,T}$ bounded by some
element of $\WW$.  In \cite{Pierre1} M. Pierre introduced the  so-called associated precise $m_1$-version $\tilde{u}$ of $u$. By the
definition, $\tilde{u}$ is the unique quasi-u.s.c. function such
that
\[
\{v\in\WW_{0,T}+\PP^2:\tilde v\ge u \mbox{
q.e.}\}=\{v\in\WW_{0,T}+\PP^2:\tilde v\ge \tilde u \mbox{ q.e.}\}
\]
In general, it is not true that $\tilde u=u$ $m_1$-a.e. However,
if $u$ is quasi-c\`adl\`ag, then
\begin{equation}
\label{eq3.12} \hat u\le\tilde u\quad\mbox{q.e.}
\end{equation}
Indeed, by \cite[Proposition IV-I]{Pierre1} there exists a
sequence $\{u_n\}\subset \WW_{0,T}$ such that
$\inf_{n\ge1}u_n=\tilde u$ q.e. Since $u\le u_n$ q.e. and $u$ is
quasi-c\`adl\`ag, we have
\[
u(\bX_t)\le u_n(\bX_t),\quad t\in [0,\zeta_\upsilon),
\]
which implies
\[
\hat{u}(\bX_{t-})=u(\bX)_{t-}\le u_n(\bX)_{t-}
=u_n(\bX_{t-}),\quad t\in (0,\zeta_\upsilon).
\]
Taking infimum, we get
\[
\hat{u}(\bX_{t-})\le \tilde u(\bX_{t-}),\quad t\in
(0,\zeta_\upsilon).
\]
By the above and \cite[Proposition 11.1]{GS} we get
(\ref{eq3.12}). Set
\begin{equation}
\label{eq.eg.eg}u_i(t,x)=u_i(t)= \left\{
\begin{array}{ll}1, &t\in I_i ,
\smallskip\\
 0, & t\in [0,T]\setminus I_i,
\end{array}
\right.
\end{equation}
where $I_1=[0,T]\cap \mathbb Q$, $I_2=[\frac T2, T]$, $I_3=(\frac T2, T]$.
Then $\tilde u_1\equiv 1$, so $\tilde u_1$ is not an $m_1$-version of $u_1$.
Moreover, $\tilde u_2=u_2$ and $\hat u_2=u_3$ ($u_3$ is not quasi-u.s.c.), so in this case $\hat u_2<\tilde u_2$ on the set $\{\frac T2\}\times E$. Since Cap($\{t\}\times B)=0$ if and only if $m(B)=0$, it follows that in general,
(\ref{eq3.12}) cannot be replaced by  ``$\hat u=\tilde u$ q.e."

\section{Reflected BSDEs}

In what follows  $(\Omega,\FF,P)$ is a probability space,
$\mathbb{F}=\{\FF_t,t\ge 0\}$ is a filtration satisfying the usual
conditions, and  $T$ is an arbitrary, but fixed bounded
$\BF$-stopping time. By $\mathcal{T}$ we denote the set of all
$\BF$-stopping times with values in $[0,T]$. For $\sigma,\gamma\in
[0,T]$ such that $\sigma\le \gamma$ we denote by
$\mathcal{T}_{\gamma}$ (resp. $\mathcal{T}_{\sigma,\gamma}$) the
set of all $\tau\in\mathcal{T}$ such that $P(\tau\in
[\gamma,T])=1$ (resp. $P(\tau\in [\sigma,\gamma])=1$).

By $\MM$ (resp. $\MM_{loc}$) we denote the space of martingales
(resp. local martingales) with respect to $\BF$. $[M]$ denotes the
quadratic variation process of $M\in\MM_{loc}$. By  $\MM_0$ (resp.
$\MM^p$ ) we denote  the subspace of $\MM$ consisting of all $M$
such that $M_0=0$ (resp. $E[M]^{p/2}_T<\infty$).

By $\VV$ (resp. $\VV^{+}$) we denote the space of all
$\mathbb{F}$-progressively measurable processes (resp. increasing
processes) $V$ of finite variation such that $V_0=0$. $\VV^p$ is
the subspace of $\VV$ consisting of $V$ such that
$E|V|^p_T<\infty$, where $|V|_t$ denotes the variation of $V$ on
the interval $[0,t]$. $^{p}\VV$ is the space of all predictable
processes in $\VV$.

By $\mathcal{S}^{p}$ we denote the space of
$\mathbb{F}$-progressively measurable processes $Y$ such that
$E\sup_{t\le T}|Y_t|^{p}<\infty$. $L^p(\FF)$ is the space of
$\mathbb{F}$-progressively measurable processes $X$ such that
$E\int_{0}^T|X_t|^p\,dt<\infty$. $L^p(\FF_T)$ is the space of
$\FF_T$-measurable random variables $X$ such that $E|X|^p<\infty$.

Let $f:\Omega\times[0,T]\times\BR\rightarrow \BR$ be a measurable
function such that $f(\cdot,y)$ is $\mathbb{F}$-progressively
measurable for every $y\in\BR$, $\xi$ be a $\FF_T$-measurable
random variable and $V$ be a c\`adl\`ag process of finite
variation such that $V_0=0$.

\begin{df}
We say that a pair of processes $(Y, M)$ is a solution of the
backward stochastic differential equation with terminal
condition $\xi$ and right-hand side $f+dV$
(BSDE$(\xi, f+dV)$ for short) if
\begin{enumerate}
\item[(a)]$Y$ is an $\mathbb{F}$-adapted c\`adl\`ag
process of Doob's class (D),
$M\in\MM_{0,loc}$\,,
\item[(b)] $[0,T]\ni t\rightarrow f(t,Y_{t})\in L^{1}(0,T)$ and
\[
Y_{t}=\xi+\int_{t}^{T}f(r, Y_{r})\,dr+ \int_{t}^{T}dV_r-
\int_{t}^{T}dM_{r}, \quad t\in[0, T],\quad P\mbox{-a.s.}
\]
\end{enumerate}
\end{df}

Let $L,U$ be two c\`adl\`ag $\mathbb{F}$-adapted processes such
that  $L_t\le U_t,\, t\in [0,T]$, and $L_T\le\xi\le U_T$.

\begin{df}
We say that a triple of processes $(Y, M, K)$ is a solution of the
reflected backward stochastic differential equation with terminal
condition $\xi$, right-hand side $f+dV$ and lower barrier $L$
(${\rm \underline{R}}$BSDE$(\xi, f+dV, L)$ for short) if
\begin{enumerate}
\item[(a)]$Y$ is an $\mathbb{F}$-adapted c\`adl\`ag
process of Doob's class (D),
$M\in\MM_{0,loc}$, $K\in{ } ^p\mathcal{V}^{+}$,
\item[(b)]$Y_{t}\ge L_{t}$, $t\in[0,T]$,
\item[(c)]
$ \int_{0}^{T}(Y_{t-}-L_{t-})\,dK_{t}=0$, $P$-a.s.,
\item[(d)] $[0,T]\ni t\rightarrow f(t,Y_{t})\in L^{1}(0,T)$ and
\[
Y_{t}=\xi+\int_{t}^{T}f(r, Y_{r})\,dr+
\int_{t}^{T}dV_r+\int_{t}^{T}dK_{r}- \int_{t}^{T}dM_{r}, \quad
t\in[0, T],\quad P\mbox{-a.s.}
\]
\end{enumerate}
\end{df}

\begin{df}
We say that a triple of processes $(Y,M,K)$ is a solution of the
reflected backward stochastic differential equation with terminal
condition $\xi$, right-hand side $f+dV$ and upper barrier $L$
(${\rm \overline{R}}$BSDE$(\xi, f+dV, U)$ for short) if the triple
$(-Y,-M,K)$ is a solution of ${\rm \underline{R}}$BSDE$(-\xi,
\tilde{f}-dV, -L)$, where $\tilde{f}(t,y)=-f(t,-y)$.
\end{df}

\begin{df}
We say that a triple of processes $(Y, M, R)$ is a solution of the
reflected BSDE with terminal condition $\xi$, right-hand side
$f+dV$, lower barrier $L$ and upper barrier $U$ (RBSDE$(\xi, f+dV,
L, U)$ for short) if
\begin{enumerate}
\item[(a)] $Y$ is an $\mathbb{F}$-adapted c\`adl\`ag process  of class
(D), $M\in\mathcal{M}_{0, loc}$, $R\in{ }^p\mathcal{V}$,
\item[(b)]$L_t\le Y_t\le U_t$, $t\in[0,T]$, $P$-a.s.
\item[(c)]
$\int_0^T(Y_{t-}-L_{t-})\,dR_t^+=\int_0^T(U_{t-}-Y_{t-})\,dR_t^-=0$,
$P$-a.s.
\item[(d)]$[0,T]\ni t\mapsto f(t, Y_t)\in L^1(0, T)$ and
\[
Y_t=\xi+\int_t^Tf(r, Y_r)\,dr+\int_t^TdV_r+\int_t^TdR_r
-\int_t^TdM_r, \quad t\in[0,T],\quad P\mbox{-a.s.}
\]
\end{enumerate}
\end{df}

\begin{uw}
\label{uw5.1} Observe that if $(Y,M,K)$ is a solution of   ${\rm
\underline{R}}$BSDE$(\xi, f+dV, L)$ and $A\in\, ^p\VV^+$ has the
property that $dA\le dK$, then the triple $(Y,M,K-A)$ is a
solution of ${\rm \underline{R}}$BSDE$(\xi, f+dV+dA, L)$.
\end{uw}

In the sequel we will need the following lemma.

\begin{lm}
\label{lm5.1} Assume that
\begin{enumerate}[{\rm (i)}]
\item there is $\mu\in\BR$ such that for a.e.
$t\in[0, T]$ and all $y,y'\in\BR$,
\[
(f(t, y)-f(t, y'))(y-y')\le\mu|y-y'|^{2},
\]
\item $[0, T]\ni t \mapsto f(t,y)\in L^1(0,T)$ for every $y\in\BR$,
\item $\BR\ni y\mapsto f(t, y)$ is
continuous for a.e. $t\in[0, T]$,
\item $\xi\in L^1(\FF_T)$, $V\in\mathcal{V}^1$, $f(\cdot, 0)\in
L^1(\FF)$.
\end{enumerate}
Let $f_n=f \vee (-n)$, and
let $(Y^n,M^n,R^n)$ be a solution  of \mbox{\rm
RBSDE}$(\xi,f_n+dV,L,U)$. Then
\[
Y^n_t\searrow Y_t,\quad M^n_t\rightarrow M_t,\quad t\in[0,T],
\quad dR^{+,n}\nearrow dR^+,\quad dR^{-,n}\searrow dR^-,
\]
where $(Y,M,R)$ is a solution of \mbox{\rm RBSDE}$(\xi,f+dV,L,U)$.
\end{lm}
\begin{dow}
By \cite[Proposition 2.14, Proposition 3.1, Theorem 3.3]{K:SPA2},
$Y^n\ge Y^{n+1}\ge Y$, $dR^{+,n}\le dR^{+,n+1}\le dR^+$,
$dR^{-,n}\ge dR^{-,n+1}$ for $n\ge 1$. Set
$\bar{Y}_t=\lim_{n\rightarrow \infty} Y^n_t,\, K_t
=\lim_{n\rightarrow \infty}R^{+,n}_t,\, A_t=\lim_{n\rightarrow
\infty}R^{-,n}_t,\, \bar{R}_t=K_t-A_t,\, t \in [0,T]$. Without
loss of generality we may assume that $\mu\le 0$ in (i). Then
\[
f(r,Y^1_r)\le f_n(r,Y^n_r)\le f_1(r,Y_r),\quad r\in [0,T].
\]
From this we conclude that the sequence $\{M^n\}$ is locally
uniformly integrable. Hence $\bar{M}_t:=\lim_{n\rightarrow \infty}
M^n_t$, $t\in [0,T]$, is a local martingale. We shall show that
the triple $(\bar{Y},\bar{M},\bar{R})$ is a solution of
RBSDE$(\xi,f+dV,L,U)$. It is clear that $\bar{R}$ is a c\`adl\`ag
process of finite variation. Moreover,  by (ii) and (iii),
\[
\bar{Y}_t=\xi+\int_t^Tf(r,\bar{Y}_r)\,dr+\int_t^T\,d\bar{R}_r
+\int_t^T\,dV_r-\int_t^T\,d\bar{M}_r,\quad t\in [0,T].
\]
It is also clear that $\bar{Y}$ is of class (D) and $L\le
\bar{Y}\le U$. By the Hahn-Saks theorem,
\begin{equation}
\label{eq4.2}
\int_0^T(\bar{Y}_t-L_t)\,dK^c_t=\lim_{n\rightarrow\infty}
\int_0^T(\bar{Y}_t-L_t)\,dR^{+,n,c}_t=0.
\end{equation}
Assume that $\Delta K_t>0$. Then there exists $n_0$ such that
$\Delta R^{+,n}_t>0$ for $n\ge n_0$. Since
$\int^T_0(Y^n_{t-}-L_{t-})\,dR^{+,n}_t=0$ and $Y^n_t\ge L_t$, it
follows that  $Y^n_{t-}=L_{t-}$ for $n\ge n_0$. Consequently,
$\bar{Y}_{t-}\le Y^{n}_{t-}=L_{t-}$, which implies that
$\bar{Y}_{t-}=L_{t-}$. We have proved that
$\sum_{t}(\bar{Y}_{t-}-L_{t-})\Delta{K}_t=0 $, which when combined
with (\ref{eq4.2}) shows that
\[
\int_0^T(\bar{Y}_{t-}-L_{t-})\,dK_t=0.
\]
Also observe that
\begin{equation}
\label{eq.uw.uw}
\int_0^T(U_t-\bar{Y}_t)\,dA_t\le \int_0^T (U_t-Y^n_t)\,dR^{n,-}_t=0.
\end{equation}
Thus the triple $(\bar{Y},\bar{M},\bar R)$ is a solution of
RBSDE$(\xi,f+dV,L,U)$. From this and the  minimality of the Jordan
decomposition of the measure $R$ it follows that $R^{+}=K,\,
R^{-}=A$.
\end{dow}

\section{PDEs with one reflecting barrier}
\label{sec5}

In this section $T>0$ is a real number, $\varphi:E\rightarrow\BR$,
$f:\bar{E}_{0,T}\times\BR\rightarrow\BR$ are measurable functions,
and $h:\bar{E}_{0,T}\rightarrow\BR$ is a quasi-c\`adl\`ag function
such that $\hat h(T,\cdot)\le\varphi$.

\begin{df}
We say that a quasi-c\`adl\`ag function $u$ on $\bar E_{0,T}$ is a
solution of the Cauchy problem
\[
-\frac{\partial u}{\partial t}-L_t u= f(\cdot,u)+\mu,\quad
u(T,\cdot)=\varphi
\]
(PDE$(\varphi,f+d\mu)$ for short)  if for q.e. $z\in E_{0,T}$,
\[
u(z)=E_z\varphi(\mathbf{X}_{\zeta_\upsilon})
+E_z\int_0^{\zeta_\upsilon}f(\mathbf{X}_r,u(\mathbf{X}_r))\,dr
+E_z\int_0^{\zeta_\upsilon}\,
dA^\mu_r.
\]
\end{df}

\begin{df}
We say that a pair $(u,\nu)$, where
$u:\bar{E}_{0,T}\rightarrow\BR$ is a quasi-c\`adl\`ag function and
$\nu$ is a positive smooth measure on $E_{0,T}$, is a solution of
problem (\ref{eq1.1}), i.e. obstacle problem with data
$\varphi,f,\mu$ and lower barrier $h$ (we denote it by ${\rm
\underline{O}P}(\varphi,f+d\mu,h)$), if
\begin{enumerate}[(a)]
\item (\ref{eq1.3}) is satisfied for some local martingale
$M$ (with $M_0=0$),
\item $u\ge h$ q.e.,
\item $\int_{E_{0,T}}(\hat{u}-\hat{h})\,d\nu=0$.
\end{enumerate}
\end{df}

\begin{df}
We say that a pair $(u,\nu)$ is a solution of the obstacle problem
with data $\varphi,f,\mu$ and upper barrier  $h$ (we denote it by
${\rm \overline{O}P}(\varphi,f+d\mu,h)$), if $(-u,\nu)$ is a
solution of ${\rm \underline{O}P}(-\varphi,\tilde{f}-d\mu,-h)$,
where $\tilde{f}(t,x,y)=-f(t,x,-y)$.
\end{df}

We will need the following duality condition considered in
\cite{K:JFA}.

\begin{enumerate}
\item[$(\Delta)$]  For some
$\alpha\ge 0$ there exists a nest $\{F_n\}$ on $E_{0,T}$ such that
for  every $n\ge 1$ there is a non-negative $\eta_n \in
L^2(E_{0,T};m_1)$ such that $\eta_n>0$ $m_1$-a.e. on $F_n$ and
$\hat{G}_\alpha^{0,T}\eta_n$ is bounded, where
$\hat{G}_\alpha^{0,T}$ is the adjoint operator to
$G_\alpha^{0,T}$.
\end{enumerate}

Note that  $(\Delta)$ is satisfied  if for some $\gamma\ge0$ the
form $\EE_{\gamma}$ has the dual Markov property (see \cite[Remark
3.9]{K:JFA}).

Following
\cite{K:JFA} we say that $f:\bar{E}_{0,T}\rightarrow\BR$ is
quasi-integrable if $f\in \BB(E_{0,T})$ and
$P_{z}(\int_{0}^{\zeta_{\tau}}|f|(\bX_{r})\,dr<\infty)=1$ for q.e.
$z\in E_{0,T}$\,.

The set of all quasi-integrable functions will be denoted by
$qL^{1}(E_{0, T};m_1)$. Note that by \cite[Remark 5.1]{K:JFA}, under condition
($\Delta$), if $f$ satisfies the condition
\begin{equation}
\label{eq6.1} \forall\, \varepsilon>0\,\, \exists\,
F_{\varepsilon}\subset E_{0,T}, \,\,
F_{\varepsilon}\mbox{-closed},\,\, \mbox{Cap}(E_{0,T}\setminus
F_{\varepsilon})<\varepsilon,\,\, \mathbf{1}_{F_{\varepsilon}}f\in
L^1(E_{0,T};m_1),
\end{equation}
then $f$ is quasi-integrable. Also note that by \cite[Proposition 3.8]{K:JFA},
under condition ($\Delta$),
\[
\MM_1\subset \mathbb{M}.
\]

We say that a function $u:\bar{E}_{0,T}\rightarrow\BR$ is of class
(D) if process $u(\bX)$ is of class (D), i.e. family $\{u(\bX_\tau),\, \tau\le \zeta_\upsilon,\, \tau-{\rm stopping\, time}\}$ is uniformly integrable under measure $P_z$ for q.e. $z\in E_{0,T}$.

\begin{uw}
\label{uw6.1} If $E_z|u(\bX_{\zeta_\upsilon})|<\infty$ for q.e. $z\in E_{0,T}$ (this holds for example if $u(T,\cdot)\in L^1(E;m)$, see \cite[Proposition 3.8]{K:JFA}) and there exists a  potential $v$ on $\bar E_{0,T}$ such that $|u|\le v$ q.e. on $E_{0,T}$ then $u$  is of class (D). Indeed, let $v=R^{0,T}\beta$ for some positive  $\beta\in\mathbb M_{T}$. By \cite[Proposition 3.4]{K:JFA}
\[
v(\mathbf{X}_t)=v(\mathbf{X}_0)-\int_0^t\,dA^\beta_r
+\int_0^t\,dM_r,\quad t\in [0,\zeta_\upsilon]
\]
for some martingale $M$. Of course, the process
$v(\mathbf{X})$ is of class (D), and $|u|(\bX_t)\le v(\bX_t),\,
t\in [0,\zeta_\upsilon)$. Since
$E_z|u|(\bX_{\zeta_\upsilon})<\infty$ for q.e. $z\in E_{0,T}$, it
follows that  $u(\bX)$ is of class (D), too.
\end{uw}

Let us consider the following assumptions.
 \begin{enumerate}[(H1)]
\item $\mu,f(\cdot,0)\cdot m_1 \in \mathbb{M},\, E_z|\varphi(\bX_{\zeta_\upsilon})|<\infty$ for q.e. $z\in E_{0,T}$,
\item there exists $\lambda\in\BR$ such that for all
$y,y'\in\BR$ and
$(t,x)\in E_{0,T}$,
\[
(f(t,x,y)-f(t,x,y'))(y-y')\le\lambda|y-y'|^2,
\]
\item the mapping $y\mapsto f(t,x,y)$ is continuous for
every $(t,x)\in E_{0,T}$,
\item $f(\cdot,y)\in qL^1(E_{0,T};m_1)$ for every $y\in\BR,$
\item $h^+$ is of class (D).
\end{enumerate}

\begin{uw}
It is an elementary check (see \cite[(3.7)]{K:JFA}) that $E_z|\varphi(\bX_{\zeta_\upsilon})|=E_zA^{|\beta|}_{\zeta_\upsilon}$ (on $E_{0,T}$), where $\beta=\delta_{\{T\}}\otimes (\varphi\cdot m)$, so by \cite[Proposition 3.8]{K:JFA}
under condition $(\Delta)$ if $\varphi\in L^1(E;m)$ then $E_z|\varphi(\bX_{\zeta_\upsilon})|<\infty$ for q.e. $z\in E_{0,T}$.
\end{uw}

\begin{tw}
\label{tw6.1} Assume {\rm (H2)}. Then there exists at most one
solution  of ${\rm\underline{O}P}(\varphi,f+d\mu,h)$.
\end{tw}
\begin{dow}
Let $(u_1,\nu_1),\, (u_2,\nu_2)$ be two solutions of
${\rm\underline{O}}$P$(\varphi,f+d\mu,h)$. Set $u=u_1-u_2$,
$\nu=\nu_1-\nu_2$  By the Tanaka-Meyer formula, for any stopping
time $\tau$  such that $\tau\le\zeta_\upsilon$ we have
\begin{align*}
|u(\mathbf{X}_t)|&\le |u(\mathbf{X}_\tau)|
+\int_t^\tau(f(r,u_1(\mathbf{X}_r))
-f(r,u_2(\mathbf{X}_r))){\rm sgn}(u)(\mathbf{X}_r)\,dr\\
&\quad +\int_t^\tau{\rm sgn}(u)(\bX_{r-})\,
dA^\nu_r-\int_t^\tau{\rm sgn}(u)(\bX_{r-})\, dM_r,\quad
t\in[0,\tau]
\end{align*}
for q.e. $z\in E_{0,T}$. Observe that by the minimality condition
and Proposition \ref{stw4.1},
\begin{align}
\label{eq5.4} \int_t^\tau{\rm sgn}(u)(\bX_{r-})\, dA^\nu_r
&=\int_t^\tau\mathbf{1}_{\{u\neq 0\}}\frac{1}{|u|} u(\bX_{r-})\,
dA^\nu_r\nonumber\\
&\le \int_t^\tau\mathbf{1}_{\{u\neq 0\}}\frac{1}{|u|}
(u_1-h)(\bX_{r-})\, dA^{\nu_1}_r\nonumber \\
&= \int_t^\tau\mathbf{1}_{\{u\neq 0\}}
\frac{1}{|u|}(\hat{u}_1(\bX_r)-\hat{h}(\bX_r))\,dA^{\nu_1}_r=0.
\end{align}
Let $\{\tau_k\}$ be a chain (i.e. an increasing sequence of stopping times with the property
$P_z(\liminf_{k\rightarrow\infty}\{\tau_k=T\})=1$, q.e.) such that
$M^{\tau_k}$ is a martingale for each $k\ge1$. Then by (H2) and
(\ref{eq5.4}),
\begin{equation}
\label{eq5.2} E_z|u(\bX_t)|\le E_z|u(\bX_{\tau_k})|+\lambda
E_z\int_t^{\tau_k}|u(\bX_r)|\,dr.
\end{equation}
Letting $k\rightarrow\infty$ in (\ref{eq5.2}) and using the fact
that $u$ is of class (D) we get
\[
E_z|u(\bX_t)|\le\lambda E_z\int_t^{\zeta_\upsilon}|u(\bX_r)|\,dr.
\]
Applying Gronwall's lemma shows that $E_z|u|(\bX_t)|=0$, $t\in
[0,\zeta_\upsilon]$ for  q.e. $z\in E_{0,T}$. This implies that
$u=0$ q.e. From this, (\ref{eq1.3}) and the uniqueness of the
Doob-Meyer decomposition we conclude that $A^\nu=0$, which forces
$\nu=0$.
\end{dow}

\begin{tw}
\label{tw6.2} Assume {\rm (H1)--(H5)}. Then there  exists a
solution of ${\rm\underline{O}P}(\varphi,f+d\mu,h)$.
\end{tw}
\begin{dow}
By \cite[Theorem 5.8]{K:JFA} there exists a unique solution $u_n$ of
(\ref{eq1.8}). Set
\[
f_n(t,x,y)=f(t,x,y)+n(y-h(t,x))^-.
\]
By the definition of a solution of (\ref{eq1.8}) and \cite[Proposition 3.4]{K:JFA} there exists a
martingale $M^n$ such that the pair $(u_n(\bX),M^n)$ is a solution of
BSDE$(\varphi(\bX_{\zeta_\upsilon}),f_n(\mathbf{X},\cdot)+dA^\mu)$
under the measure $P_z$ for q.e. $z\in E_{0,T}$.  By \cite[Theorem
4.1]{K:SPA2} (see also \cite{T}) there exists a solution $(Y^z,M^z,K^z)$ of ${\rm
\underline{R}}$BSDE$(\varphi(\mathbf{X}),f(\bX,\cdot)+dA^\mu,h(\bX))$
under the measure $P_z$ for q.e. $z\in E_{0,T}$, and
\begin{equation}
\label{eq5.3} u_n(\bX_t)\nearrow Y^z_t,\quad t\in
[0,\zeta_\upsilon],\quad P_z\mbox{-a.s.}
\end{equation}
From (\ref{eq5.3}) it follows that  $u_n\le u_{n+1}$ q.e. (since the exceptional sets coincide  with the sets of  zero capacity) for
$n\ge1$ and
\[
Y^z_t=u(\bX_t),\quad t\in [0,\zeta_\upsilon],\quad P_z\mbox{-a.s.}
\]
for q.e. $z\in E_{0,T}$, where $u:=\sup_{n\ge 1} u_n$.  This
implies that $u$ is quasi-c\`adl\`ag. What is left to show that
there exists a smooth measure $\nu$ such that $K^z=A^\nu$ for q.e.
$z\in E_{0,T}$. Since the pointwise limit of additive functionals
is an additive functional, we may assume by Lemma \ref{lm5.1} and
\cite[Theorem 2.13]{K:SPA2} that $EK^z_{\zeta_\upsilon}<\infty$ for q.e. $z\in
E_{0,T}$. By \cite[Proposition 4.3]{K:SPA2} and Theorem \ref{tw4.1}, for every
predictable stopping time $\tau$,
\[
\Delta K^z_\tau=(u(\bX_\tau)-\hat{h}(\bX_\tau)+\Delta
A^\mu_\tau)^{-}.
\]
Hence
\[
J_t=\sum_{s\le t}\Delta K^z_t
\]
is a PNAF (without jump in $\zeta_\upsilon$ since $\hat h(T,\cdot)\le \varphi$), which implies that there exists $\beta\in\mathbb{M}$
such that $J=A^\beta$ (by the Revuz duality). Set
\[
C^z_t=K^z_t-A^\beta_t, \quad t\in
[0,\zeta_\upsilon].
\]
It is clear that the process $C^z$ is continuous. By Remark
\ref{uw5.1}, the triple $(Y^z,M^z,C^z)$ is a solution of ${\rm
\underline{R}}$BSDE$(\varphi(\mathbf{X}),f(\bX,\cdot)+dA^\mu+dA^\beta,
h(\bX))$. By \cite[Theorem 5.8]{K:JFA} there exists a solution $v_n$ of the
equation
\[
-\frac{\partial v_n}{\partial t}-L_t v_n=f_n(\cdot,v_n)+\mu+\beta,\quad v_n(T,\cdot)=\varphi.
\]
Therefore, by the definition of solution and \cite[Proposition 3.4]{K:JFA} there exists a martingale $N^n$ such that the pair $(v_n(\bX),
N^n)$ is a solution of
BSDE$(\varphi(\mathbf{X}),f_n(\bX,\cdot)+dA^\mu+dA^\beta)$ under
the measure $P_z$ for q.e. $z\in E_{0,T}$. Let
$C^n_t=\int_0^tn(v_n(\bX_r)-h(\bX_r))^-\,dr$,
$t\in[0,\zeta_\upsilon]$. By \cite[Theorem 2.13]{K:SPA2} the sequence $\{C^n\}$
converges uniformly on $[0,\zeta_\upsilon]$ in probability $P_z$
for q.e. $z\in E_{0,T}$. Since $C^n$ is a PCAF for each $n\ge1$,
the process $C$ defined by
\[
C_t=\lim_{n\rightarrow \infty} C^n_t,\quad t\in
[0,\zeta_\upsilon],
\]
is a PCAF.  Therefore there exists a measure $\alpha\in\mathbb{M}$
such that $C=A^\alpha$. It is clear that $C^z=A^\alpha$ for q.e.
$z\in E_{0,T}$. Finally, $K^z=A^\nu$ for q.e. $z\in E_{0,T}$, where
$\nu=\alpha+\beta$.
\end{dow}

\begin{stw}
\label{stw.cr}
Assume that $\varphi_i,f_i,\mu_i$, $i=1,2$, satisfy \mbox{\rm(H1)--(H4)} and $h_1,h_2$
are quasi-c\`adl\`ag. Let $(u_i,\nu_i)$ be a solution to ${\rm\underline{O}P}(\varphi_i,f_i+d\mu_i,h_i)$.
Assume that $\varphi_1\le\varphi_2$ $m$-a.e., $f_1(\cdot,y)\le f_2(\cdot,y)$ $m_1$-a.e. for every $y\in\BR$, $d\mu_1\le d\mu_2$ and $h_1\le h_2$ q.e. Then
\[
u_1\le u_2\quad\mbox{q.e. on  } E_{0,T},
\]
and if $h_1=h_2$ q.e., then $d\nu_1\ge d\nu_2$.
\end{stw}
\begin{dow}
By the proof of Theorem \ref{tw6.2}, $u_{i,n}\nearrow u_i$ q.e., where $u_{i,n}$ is a solution to
PDE$(\varphi_i,f_{n,i}+d\mu_i)$ with
\[
f_{i,n}(t,x,y)=f_i(t,x,y)+n(y-h_i(t,x))^-.
\]
By \cite[Corollary 5.9]{K:JFA}, $u_{1,n}\le u_{2,n}$ q.e. Hence $u_1\le u_2$ q.e. As for the second assertion of the theorem, by Lemma \ref{lm5.1} we may assume that $\nu_1,\nu_2\in \mathbb M$. By the proof of Theorem \ref{tw6.2},
\[
A^{\nu_i}_t=A^{\alpha_i}_t+A^{\beta^i}_t,\quad t\in [0,\zeta_\upsilon],\quad P_z\mbox{-a.s.}
\]
for q.e. $z\in E_{0,T}$ and some positive $\beta_i\in \mathbb M$, $\alpha_i\in\mathbb M_c$ such that $A^{\beta_i}$
is purely jumping and
\[
\Delta A^{\beta_i}_t=(u_i(\bX_t)-\hat{h}(\bX_t)+\Delta
A^{\mu_i}_t)^{-}.
\]
We already know that $u_1\le u_2$ q.e. Moreover, by the assumptions and Revuz duality, $dA^{\mu_1}\le dA^{\mu_2}$. Therefore
$dA^{\beta_1}\ge dA^{\beta_2}$. Furthermore, by the proof of Theorem \ref{tw6.2},
\[
A^{\alpha_{i,n}}_t\rightarrow A^{\alpha_i}_t,\quad t\in [0,\zeta_\upsilon],\quad P_z\mbox{-a.s.}
\]
for q.e. $z\in E_{0,T}$, where $\alpha_{i,n}=n(v_{i,n}-h)^-\cdot m_1$ and $v_{i,n}$ is a solution to
PDE$(\varphi_i, f_{n}+d\mu_i+d\beta_{i})$. By \cite[Corollary 5.9]{K:JFA}, $v_{1,n}\le v_{2,n}$, which implies that
$\alpha_{1,n}\ge \alpha_{2,n}$, $n\ge 1$. Consequently, $dA^{\alpha_{1,n}}\ge dA^{\alpha_{2,n}}$, $n\ge1$.
Hence $dA^{\alpha_{1}}\ge dA^{\alpha_{2}}$, so $dA^{\nu_{1}}\ge dA^{\nu_{2}}$, which implies that $d\nu_1\ge d\nu_2$.
\end{dow}

\begin{uw}
\label{uw6.2}
Let  $v=R^{0,T}\beta$ for some $\beta\in\mathbb M_T$ be such that $v\ge h$ q.e. on $E_{0,T}$.
By \cite[Proposition 3.4]{K:JFA}
\[
v(\mathbf{X}_t)=v(\mathbf{X}_0)-\int_0^t\,dA^\beta_r
+\int_0^t\,dM_r,\quad t\in [0,\zeta_\upsilon]
\]
for some martingale $M$.  Let $\gamma=-(h(T,\cdot)\cdot
m)\otimes\delta_{\{T\}}$. Observe that the function $\bar{v}=
E_\cdot h(\bX_{\zeta_\upsilon})+R^{0,T}\gamma+R^{0,T}\beta$ is
equal to $v$ on $E_{0,T}$. Moreover, $\bar{v}(X)$ satisfies the
same equation as $v(X)$, but on $[0,\zeta_\upsilon)$, and
$\bar{v}(T,\cdot)\ge h(T,\cdot)$, from which it follows that
$\bar{v}(\bX_t)\ge h(\bX_t)$, $t\in [0,\zeta_\upsilon]$.
\end{uw}

We denote by $\MM_{1,T}$ the set of all finite Borel measures on $\bar E_{0,T}$.

\begin{stw}
\label{stw6.1} Let the hypotheses {\rm(H1)--(H4)} hold. Assume that
$\mu\in\mathbb{M}$ and there exists $\beta\in \mathbb{M}_T$ such
that $v:=R^{0,T}\beta\ge h$ on $E_{0,T}$, and $f^-(\cdot,v)\cdot
m_1\in \mathbb{M}$. Let $(u,\nu)$ be a solution of
${\rm\underline{O}P}(\varphi,f+d\mu,h)$. Then $f(\cdot,u)\cdot
m_1, \nu\in\mathbb{M}$. If we assume additionally that
$\EE_\gamma$ has the dual Markov property for some $\gamma\ge 0$
and $\varphi\in L^1(E;m)$, $f(\cdot,0),\, f^-(\cdot,v)\in
L^1(E_{0,T};m_1), \mu\in\MM_1,\,  \beta\in\MM_{1,T}$, then $f(\cdot,u)\in
L^1(E_{0,T};m_1)$ and $\nu\in\MM_1$.
\end{stw}
\begin{dow}
Assume that $\mu, f^-(\cdot,v)\cdot m_1\in\mathbb{M}$ and $\beta\in\mathbb M_{T}$. By
\cite[Theorem 2.13]{K:SPA2} (see also Remark \ref{uw6.2}),
$\nu\in\mathbb{M}$. By \cite[Theorem 5.4]{K:JFA},
 $f(\cdot,u)\cdot m_1\in \mathbb{M}$, which proves the first
part of the proposition. Now, assume that $\mu,
f^-(\cdot,v)\cdot m_1\in\MM_1$, $\beta\in\MM_{1,T}$. Let $\bar v$ be a solution of  PDE$(\varphi^+,f+f^-(\cdot,v)+d\mu^++d\beta^+)$ (it exists by  \cite[Theorem 5.8]{K:JFA}). By \cite[Corollary 5.9]{K:JFA},
$\bar{v}\ge v$,  and consequently $\bar v\ge h$ q.e. on $E_{0,T}$. Observe that the pair $(\bar v, f^-(\cdot,v)\cdot m+\beta^+)$ is a solution to ${\rm\underline{O}P}(\varphi^+,f+d\mu^+,\bar v)$. Hence, by Proposition \ref{stw.cr},
$\bar v\ge u$ q.e. on $E_{0,T}$. On the other hand, by \cite[Corollary 5.9]{K:JFA}, $u_0\le u$ q.e. on $E_{0,T}$,
where $u_0$ is a solution to PDE$(\varphi,f+d\mu)$.
By \cite[Proposition 5.10]{K:JFA}, $f(\cdot,\bar{v}),\,
f(\cdot,u_0), \bar{v}, u_0\in L^1(E_{0,T}; m_1)$. Since $u_0\le
u\le \bar{v}$, thanks to (H2) we have $u, f(\cdot,u)\in L^1(E_{0,T};
m_1)$.  Let $\gamma=f^-(\cdot,v)\cdot m+\mu^++\beta^+$.
Observe that
\begin{align*}
u(z)&=E_z\varphi(\bX_{\zeta_\upsilon})
+E_z\int_0^{\zeta_\upsilon}f(\bX_r,u(\bX_r))\,dr
+E_z\int_0^{\zeta_\upsilon}\,dA^\nu_r+E_z\int_0^{\zeta_\upsilon}\,dA^\mu_r\\
&\le E_z\varphi^+(\bX_{\zeta_\upsilon})
+E_z\int_0^{\zeta_\upsilon}f(\bX_r,\bar{v}(\bX_r))\,dr
+E_z\int_0^{\zeta_\upsilon}\,dA^\gamma_r=\bar{v}(z).
\end{align*}
Hence
\begin{align*}
E_z\int_0^{\zeta_\upsilon}\,dA^\nu_r &\le E_z|\varphi|(\bX_{\tau_\upsilon})+
E_z\int_0^{\zeta_\upsilon}|f(\bX_r,u(\bX_r))|\,dr
\\&\quad+E_z\int_0^{\zeta_\upsilon}|f(\bX_r,\bar{v}(\bX_r))|\,dr
+E_z\int_0^{\zeta_\upsilon}\,dA^{|\mu+\gamma|}_r
\end{align*}
for q.e. $z\in E_{0,T}$. By the above inequality and
\cite[Proposition 3.13]{K:JFA},
\[
\|\nu\|_1\le
c(\|\varphi\|_{L^1}+\|f(\cdot,u)\|_{L^1}+\|f^-(\cdot,v)\|_{L^1}+\|\mu^+\|_1+\|\beta^+\|_1),
\]
which completes the proof.
\end{dow}

\begin{uw}
\label{uw6.l1}
Observe that under the assumptions of the second assertion of the
above proposition, $u\in L^1(E_{0,T}; m_1)$. Moreover, by
\cite[Proposition 5.10]{K:JFA},
\[
\|u\|_{L^1}+\|f(\cdot,u)\|_{L^1}\le c(\|\varphi\|_{L^1}
+\|f(\cdot,0)\|_{L^1}+\|\mu\|_1+\|\nu\|_1).
\]
\end{uw}

\section{PDEs with two reflecting barriers}
\label{sec6}

We assume as given $T,\varphi, f$ as in Section \ref{sec5}, and
quasi-c\`adl\`ag functions $h_1,h_2:\bar E_{0,T}\rightarrow\BR$
such that $\hat h_1(T,\cdot)\le\varphi\le \hat h_2(T,\cdot)$.

\begin{df}
We say that a pair $(u,\nu)$ consisting of a quasi-c\`adl\`ag
function $u:\bar E_{0,T}\rightarrow\BR$ of class (D) and a smooth
measure $\nu$ on $E_{0,T}$ is a solution of the obstacle problem
with data $\varphi,f,\mu$ and barriers $h_1,h_2$ (we denote it by
OP$(\varphi,f+d\mu,h_1,h_2))$ if
\begin{enumerate}[(a)]
\item (\ref{eq1.3}) is satisfied for some local martingale $M$ (with $M_0=0$),
\item $h_1\le u\le h_2$ q.e. on $E_{0,T}$,
\item $\int(\hat{u}-\hat{h}_1)\,d\nu^+=\int(\hat{h}_2-\hat{u})\,d\nu^-=0.$
\end{enumerate}
\end{df}

\begin{stw}
\label{stw.unq}
Let assumption  \mbox{\rm(H1)} hold. Then there exists at most one
solution of \mbox{\rm OP}$(\varphi,f+d\mu,h_1,h_2)$.
\end{stw}
\begin{dow}
The proof is analogous to the proof of Theorem \ref{tw6.1}. The
only difference  is the estimate of the integral involving
$dA^\nu$. In the present situation, the estimate is as follows:
\begin{align}
\label{eq.uqn12}
\nonumber\int_t^\tau{\rm sgn}(u)(\bX_{r-})\, dA^\nu_r&\nonumber
\le \int_t^\tau\mathbf{1}_{\{u\neq 0\}}
\frac{1}{|u|} (u_1-h_1)(\bX_{r-})\, dA^{\nu_1}_r\\
&\nonumber\quad
+ \int_t^\tau\mathbf{1}_{\{u\neq 0\}}
\frac{1}{|u|} (h_2-u_2)(\bX_{r-})\, dA^{\nu_2}_r\\
&=\int_t^\tau \mathbf{1}_{\{u\neq 0\}}
\frac{1}{|u|}(\hat{u}_1(\bX_r)-\hat{h}_1(\bX_r))\,dA^{\nu_1}_r\nonumber\\
&\quad
+ \int_t^\tau \mathbf{1}_{\{u\neq 0\}}
\frac{1}{|u|}(\hat{h}_2(\bX_r)-\hat{u}_2(\bX_r))\,dA^{\nu_2}_r =0.
\end{align}
\end{dow}

Consider the following hypothesis:

\begin{enumerate}
\item[(H6)] $h_1, h_2$ are quasi-c\`adl\`ag functions of class (D) such that
$\hat h_1(T,\cdot)\le \varphi\le \hat h_2(T,\cdot)$ and $h_1<h_2$,
$\hat{h}_1<\hat{h}_2$ q.e. on $E_{0,T}$, or there exists
$\beta\in\mathbb{M}_T$ such that $h_1\le R^{0,T}\beta\le h_2$ q.e.  on
$E_{0,T}$.
\end{enumerate}

\begin{uw}
\label{uw7.1} If (H1)--(H4) (H6) are satisfied then the assertion of \cite[Theorem 4.2]{K:SPA2}
holds true under measure $P_z$ for q.e. $z\in E_{0,T}$. Indeed, it is clear that (H1)--(H4) of \cite{K:SPA2} are satisfied under measure $P_z$ for q.e. $z\in E_{0,T}$.
If $h_1<h_2$ and $\hat{h}_1<\hat{h}_2$ q.e. on
$E_{0,T}$, then (since the exceptional sets coincide  with the sets of  zero capacity) for q.e. $z\in E_{0,T}$
\begin{equation}
\label{eq6.2} h_1(\bX_t)< h_2(\bX_t),\quad
\hat{h}_1(\bX_t)<\hat{h}_2(\bX_t), \quad t\in (0,\zeta_\upsilon),\quad P_z\mbox{-a.s.}
\end{equation}
Since $\bX$ has no predictable jumps (as a Hunt process), it follows from
(\ref{eq6.2}) and Theorem \ref{tw4.1} that
$h_1(\bX)_{\tau-}<h_2(\bX)_{\tau-}$ for every predictable stopping
time $\tau$ with values in $(0,\zeta_\upsilon]$. By this and
\cite{T} the assertion of \cite[Theorem 4.2]{K:SPA2} is satisfied under measure $P_z$ for q.e. $z\in E_{0,T}$. If
$h_1\le R^{0,T}\beta\le h_2$ on $E_{0,T}$ for some
$\beta\in\mathbb{M}_T$, then hypothesis (H7) from \cite{K:SPA2} (with $L=h_1(\bX),\, U=h_2(\bX),\, X=v(\bX)$) is
satisfied by  Remark \ref{uw6.2} under measure $P_z$ for q.e. $z\in E_{0,T}$, so again the assertion of \cite[Theorem 4.2]{K:SPA2} holds true.
\end{uw}

\begin{tw}
\label{tw7.1} Assume that {\rm (H1)--(H4), (H6)} are satisfied.
Then there  exists a unique solution of {\rm
OP}$(\varphi,f+d\mu,h_1,h_2)$.
\end{tw}

\begin{dow}
By Remark \ref{uw7.1}, for q.e. $z\in E_{0,T}$  there
exists a solution $(Y^z,M^z,R^z)$ of
RBSDE$(\varphi(\bX_{\tau_\upsilon}), f(\bX,\cdot)+dA^\mu,
h_1(\bX),h_2(\bX))$ under the measure $P_z$. To prove the
existence of a solution, it suffices to show that there exists a
function $u$ and a smooth measure $\nu$ such that $Y^z=u(\bX),\,
R^z=A^\nu$ for q.e. $z\in E_{0,T}$, because then the pair
$(u,\nu)$ will be a solution of {\rm
OP}$(\varphi,f+d\mu,h_1,h_2)$. By \cite[Theorem 4.1]{K:SPA2}, for
every $n\ge 1$ there exists  a solution
$(\bar{Y}^{n,z},\bar{M}^{n,z},\bar{A}^{n,z})$ of
${\rm\overline{R}BSDE}(\varphi(\bX_{\zeta_\upsilon}),f_n(\bX,\cdot)
+dA^\mu,h_2(\bX))$ with $f_n$ defined by
\[
f_n(z,y)=f(z,y)+n(y-h_1(z))^{-}.
\]
By Theorem \ref{tw6.2},
\begin{equation}
\label{eq7.1}
\bar{Y}^{n,z}=\bar{u}_n(\bX),\quad
\bar{A}^{n,z}=A^{\bar{\gamma}_n},
\end{equation}
where  $(\bar{u}_n,\bar{\gamma}_n)$ is a solution of ${\rm
\overline{O}P}(\varphi,f_n+d\mu,h_2)$. By  (\ref{eq7.1}),
$\bar{M}^{n,z}$ in fact does not depend on $z$.  By \cite[Theorem 4.2]{K:SPA2},
\[
\bar{Y}^{n,z}_t\nearrow Y^z_t,\quad t\in [0,\zeta_\upsilon].
\]
Since the exceptional sets coincide  with the sets of  zero capacity, this implies that  $\bar u_n\le\bar u_{n+1}$ q.e. for $n\ge1$, and
\[
Y^z_t=u(\bX_t),\quad t\in [0,\zeta_\upsilon]
\]
for q.e. $z\in E_{0,T}$, where $u:=\sup_{n\ge 1}\bar{u}_n$.  By
\cite[Theorem 4.2]{K:SPA2}, $dA^{\bar{\gamma}_n}\le
dA^{\bar{\gamma}^{n+1}}$, so $A$ defined as
$A_t=\lim_{n\rightarrow \infty}A^{\bar{\gamma}_n}_t$, $t\ge0$, is
a PNAF. Therefore there exists a positive smooth measure $\gamma$
such that $A=A^\gamma$. By \cite[Theorem 4.2]{K:SPA2},
$A^\gamma=R^{z,-}$ for q.e. $z\in E_{0,T}$. By  Lemma \ref{lm5.1},
without loss of generality we may assume that
$E_z\int_0^{\zeta_\upsilon}\,d|R^z|_r<\infty$ for q.e. $z\in
E_{0,T}$. Observe that the triple $(Y^z,M^z, R^{z,+})$ is a
solution of ${\rm \underline{O}P}(\varphi,f+d\mu+d\gamma,h_1)$.
Therefore, by Theorem \ref{tw6.2}, there exists a positive smooth
measure $\alpha$ such that $R^{z,+}=A^\alpha$ for q.e. $z\in
E_{0,T}$, which completes the proof.
\end{dow}

\begin{stw}
\label{stw7.1} Let the hypotheses {\rm(H1)--(H4)} and \mbox{\rm(H6)}
with some measure $\beta\in\BM_T$ hold.  Assume that
$f(\cdot,v)\cdot m_1\in\mathbb{M}$ with $v=R^{0,T}\beta$, and that
$\mu\in\mathbb{M}$. Then $f(\cdot,u)\cdot m_1, \nu\in\mathbb{M}$.
If, in addition, $\EE_\gamma$ has the dual Markov property for
some $\gamma\ge 0$, and  $\varphi\in L^1(E;m), f(\cdot,0),
f(\cdot,v)\in L^1(E_{0,T}; m_1),\, \mu\in \MM_1$, $\beta\in \MM_{1,T}$, then
$f(\cdot,u)\in L^1(E_{0,T},m_1), \nu\in\MM_1$.
\end{stw}
\begin{dow}
If $f(\cdot,v)\cdot m_1, \mu\in\mathbb{M}$, then
$\nu\in\mathbb{M}$ by \cite[Theorem 3.3]{K:SPA2} (see also Remark
\ref{uw6.2}). Hence, by \cite[Theorem 5.4]{K:JFA}, $f(\cdot,u)\cdot m_1\in \mathbb M$. Assume that $\EE_\gamma$ has the dual Markov
property for some $\gamma\ge 0$, and that $f(\cdot,0),
f(\cdot,v)\in L^1(E_{0,T}; m_1)$, $\mu\in \MM_1$, $\beta\in\MM_{1,T}$.
Let $(\bar v,\bar \nu)$ be a solution to ${\rm
\overline{O}P}(\varphi^+,f+f^-(\cdot,v)+d\beta^++d\mu^+,h_2)$.
By Proposition \ref{stw.cr},
$\bar v\ge v$ q.e. on $E_{0,T}$ (since $(v,0)$ is a solution to ${\rm
\overline{O}P}(0,f-f(\cdot,v)+d\beta,v)$). Let $(\bar u_n,\bar\gamma_n)$ be  a solution to ${\rm\overline{O}P}(\varphi,f_n+d\mu,h_2)$ with
\[
f_n(z,y)=f(z,y)+n(y-h_1(z))^{-}.
\]
By the proof of Theorem \ref{tw7.1},
$\bar u_n\searrow u$ and $A^{\bar\gamma_n}_t\nearrow A^{\nu^-}_t,\, t\in [0,\zeta_\upsilon],\, P_z$-a.s.
for q.e. $z\in E_{0,T}$.
Since $\bar v\ge v$, we have $\bar v\ge h_1$ q.e. on $E_{0,T}$. Therefore
$f(\cdot,v)=f_n(\cdot,v)$, and in fact, $(\bar v,\bar \nu)$ is a solution to
${\rm\overline{O}P}(\varphi^+,f_n+f^-(\cdot,v)+d\beta^++d\mu^+,h_2)$.
Hence, by Proposition \ref{stw.cr},
\[
\bar\gamma_n\le \bar\nu.
\]
By the Revuz duality, $A^{\bar\gamma_n}_t\le A^{\bar\nu}_t$, which when combined with the convergence of
$\{A^{\bar\gamma_n}\}$ implies that $A^{\nu^-}_t\le A^{\bar\nu}_t$, $t\in [0,\zeta_\upsilon]$, $P_z$-a.s. for
q.e. $z\in E_{0,T}$. Consequently, $\nu^-\le \bar\nu$.
By Proposition \ref{stw6.1}, $\bar\nu\in \MM_1$. Hence, by
\cite[Proposition 3.13]{K:JFA}, $\nu^-\in \MM_1$. Since the pair
$(u,\nu^+)$ is a solution of ${\rm
\underline{O}P}(\varphi,f+d\mu-d\nu^-,h_1)$, applying Proposition
\ref{stw6.1} yields $\nu^+\in\MM_1$. This completes the proof.
\end{dow}

\begin{uw}
\label{uw7.l1}
Under  \mbox{\rm(H1)--(H4), (H6)}   and the assumptions of the
second assertion of Proposition \ref{stw7.1},
\[
\|u\|_{L^1}+\|f(\cdot,u)\|_{L^1}\le c(\|\mu\|_1+\|\varphi\|_{L^1}
+\|\nu\|_1+\|f(\cdot,0)\|_{L^1}).
\]
This follows from Proposition \ref{stw7.1} and \cite[Proposition
5.10]{K:JFA}.
\end{uw}

\begin{uw}
\label{uw.main} Let $v$ be a  difference of potentials on $\bar E_{0,T}$, i.e.
$v=R^{0,T}\beta$ for some  $\beta\in\mathbb{M}_T$. Observe that if a
pair $(\bar u,\bar \nu)$ is a solution of
OP$(\varphi,f_v+d\mu-d\beta,h_1-v,h_2-v)$ with
\[
f_v(z,y)=f(z,v+y),
\]
then  $(u,\nu)=(\bar u+v,\bar \nu)$ is a solution of
OP$(\varphi,f+d\mu,h_1,h_2)$. It follows that the assumption $f(\cdot,0)\cdot
m_1\in\mathbb{M}$ in Theorems \ref{tw6.2} and \ref{tw7.1} may be
replaced by more general assumption saying that there exists a
function $v$ which is a difference of potentials on $\bar E_{0,T}$ such that $f(\cdot,v)\cdot m_1\in\mathbb{M}$. Similarly, counterparts of the results stated
in Propositions \ref{stw6.1}, \ref{stw7.1} and Remarks
\ref{uw6.l1} and \ref{uw7.l1}  hold true under the assumption
$f(\cdot,v)\cdot m_1\in\mathbb{M}$. To get appropriate
modifications of these results, we first apply them to the pair
$(\bar u,\bar \nu)$,  and next we use the fact that
$(\bar{u},\bar{\nu})=(u-v,\nu)$.
\end{uw}

In condition (b) of the definition of the obstacle problem given at the beginning of this section we require that the solution lies q.e.  between the barriers. Below we show that if we weaken (b) and require only that this property holds a.e., then our results also apply  to measurable barriers $h_1,h_2$ satisfying the following condition: there exits  $v$ of the form $v=R^{0,T}\beta$ with  $\beta\in \mathbb M_T$ (i.e. $v$ is a difference of potentials) such that $h_1\le v\le h_2$ $m_1$-a.e. (in case of one barrier $h$ it is enough to assume that $h^+$ is of class (D)).

Before presenting  our results for measurable barriers, we  give a definition of a solution.

\begin{df}
We say that a pair $(u,\nu)$ consisting of a quasi-c\`adl\`ag
function $u:\bar E_{0,T}\rightarrow\BR$ of class (D) and a smooth
measure $\nu$ on $\bar E_{0,T}$ is a solution of the obstacle problem
with data $\varphi,f,\mu$ and measurable barriers $h_1,h_2:\bar E_{0,T}\rightarrow\BR$  if
\begin{enumerate}
\item[\rm{(a)}] (\ref{eq1.3}) is satisfied for some local martingale $M$ (with $M_0=0$),
\item[(b*)] $h_1\le u\le h_2$ $m_1$-a.e. on $E_{0,T}$,
\item[(c*)] $\int(\hat{u}-\hat{\eta}_1)\,d\nu^+=\int(\hat{\eta}_2-\hat{u})\,d\nu^-=0$ for all quasi-c\`adl\`ag
$\eta_1,\eta _2$ such that $h_1\le\eta_1\le u\le\eta_2\le h_2$ $m_1$-a.e.
\end{enumerate}
\end{df}

If the barriers are quasi-c\`adl\`ag, then the above definition agrees with
the definition given at the beginning of Section \ref{sec6}, because  for quasi-c\`adl\`ag functions condition (b*) is equivalent to (b)
and in (c*) we may take $\eta_1=h_1$ and $\eta_2=h_2$. Therefore the obstacle problem with data $\varphi,f,\mu$ and measurable barriers $h_1,h_2$ we still denote by \mbox{\rm OP}$(\varphi,f+d\mu,h_1,h_2)$.

In case of measurable barriers the proof of  uniqueness of solutions to the problem \mbox{\rm OP}$(\varphi,f+d\mu,h_1,h_2)$ goes as in the case of
quasi-c\`adl\`ag barriers, the only difference being in the fact that in (\ref{eq.uqn12}) we replace  $h_1$ by  $u_1\wedge u_2$
and $h_2$ by $u_1\vee u_2$, and then we apply (c*).

The problem of existence of a solution is more delicate.
Observe that if $(u,\nu)$ is a solution to
\mbox{\rm OP}$(\varphi,f+d\mu,h_1,h_2)$ with measurable barriers $h_1$ and $h_2$, then it is a solution to  \mbox{\rm OP}$(\varphi,f+d\mu,\eta_1,\eta_2)$ with   $\eta_1,\eta_2$  as in condition (c*). It appears that under the additional assumption on barriers (mentioned before the last definition) one can construct quasi-c\`adl\`ag $\eta_1,\eta_2$ such that if $(u,\nu)$ is a solution to
\mbox{\rm OP}$(\varphi,f+d\mu,\eta_1,\eta_2)$, then it is a solution to \mbox{\rm OP}$(\varphi,f+d\mu,h_1,h_2)$. This shows that as  long as we only require (b${}^*$), the  study of the obstacle problem with measurable barriers can be reduced to the study of
the obstacle problem with quasi-c\`adl\`ag barriers.
Finally, note that, unfortunately, there is no construction of $\eta_1,\eta_2$ depending only on $h_1$ and $h_2$ (the barriers $\eta_1,\eta_2$ depend also on $\varphi,f,\mu$).
The reason is that the class of c\`adl\`ag functions is
neither inf-stable nor sup-stable.

Let $L, U$ be  measurable adapted processes $(L\le U)$. Following  \cite{K:SPA2}
we say that a triple of processes $(Y, M, R)$ is a solution of the
reflected BSDE with terminal condition $\xi$, right-hand side
$f+dV$, lower barrier $L$ and upper barrier $U$   if
\begin{enumerate}
\item[(a)] $Y$ is an $\mathbb{F}$-adapted c\`adl\`ag process  of class
(D), $M\in\mathcal{M}_{0, loc}$, $R\in{ }^p\mathcal{V}$,
\item[(b*)]$L_t\le Y_t\le U_t$ $P$-a.s. for a.e. $t\in [0,T]$,
\item[(c*)]
$\int_0^T(Y_{t-}-H^1_{t-})\,dR_t^+=\int_0^T(H^2_{t-}-Y_{t-})\,dR_t^-=0$
$P$-a.s. for all c\`adl\`ag processes $H^1,H^2$ such that $L_t\le H^1_t\le Y_t\le H^2_t\le U_t$ $P$-a.s. for a.e. $t\in [0,T]$,
\item[(d)]$[0,T]\ni t\mapsto f(t, Y_t)\in L^1(0, T)$ and
\[
Y_t=\xi+\int_t^Tf(r, Y_r)\,dr+\int_t^TdV_r+\int_t^TdR_r
-\int_t^TdM_r, \quad t\in[0,T],\quad P\mbox{-a.s.}
\]
\end{enumerate}
It is obvious that for c\`adl\`ag barriers the above definition agrees with the definition given in Section 4.
We see that $(u,\nu)$ is a solution to  \mbox{\rm OP}$(\varphi,f+d\mu,h_1,h_2)$  if and only if $(u(\bX),A^\nu,M)$ is a solution to RBSDE$(\varphi(\bX_{\zeta_\upsilon}), f(\bX,\cdot)+dA^\mu,h_1(\bX),h_2(\bX))$ under the measure $P_z$ for q.e. $z\in E_{0,T}$ (because $m(A)=0$ if and only if $E_z\int_0^{\zeta_\upsilon}\mathbf{1}_{A}(\bX)=0$ for q.e. $z\in E_{0,T}$, q.e., and the last condition is satisfied if and only if for q.e. $z\in E_{0,T}$ we have  $\mathbf{1}_{A}(\bX_t)=0$ $P_z$-a.s. for a.e. $t\in [0,\zeta_\upsilon]$).

\begin{stw}
\label{stw.red}
Assume that \mbox{\rm(H1)--(H4)} are satisfied and $h_1,h_2$ are measurable functions such that  $h_1\le v\le h_2$ $m_1$-a.e. for some
function $v$ being a difference of potentials on $\bar E_{0,T}$.
Then there exist quasi-c\`adl\`ag functions
$\eta_1,\eta_2$ such that $h_1\le \eta_1\le v\le \eta_2\le h_2$ $m_1$-a.e., and moreover, having the property that if $(u,\nu)$ is a solution
to \mbox{\rm OP}$(\varphi,f+d\mu,\eta_1,\eta_2)$, then  $(u,\nu)$ is a solution to \mbox{\rm OP}$(\varphi,f+d\mu,h_1,h_2)$.
\end{stw}
\begin{dow}
By Remark \ref{uw7.1}, for q.e. $z\in E_{0,T}$  there exists a unique solution  $(Y^z,M^z,R^z)$ to RBSDE$(\varphi(\bX_{\zeta_\upsilon}), f(\bX,\cdot)+dA^\mu,h_1(\bX),h_2(\bX))$ under the measure $P_z$.
By \cite[Theorem 4.2]{K:SPA2},
\[
Y^{z,n}_t\rightarrow Y^z_t,\quad t\in [0,\zeta_\upsilon],\quad P_z\mbox{-a.s.}
\]
for q.e. $z\in E_{0,T}$, where $(Y^{z,n},M^{z,n})$ is a solution to BSDE$(\varphi(\bX_{\zeta_\upsilon}),f_n(\bX,\cdot)+dA^\mu)$ with
\[
f_n(z,y)=f(z,y)+n(y-h_1(z))^{-}-n(y-h_2(z))^+.
\]
By \cite[Theorem 5.8]{K:JFA}, $Y^{z,n}_t=u_n(\bX_t),\, t\in [0,\zeta_\upsilon],\, P_z$-a.s. for q.e. $z\in E_{0,T}$, where
$u_n$ is a solution to PDE$(\varphi,f_n+d\mu)$. Let $w:=\sup_{n\ge 1} u_n$. It is clear that
\[
w(\bX_t)=Y^{z}_t,\quad t\in [0,\zeta_\upsilon],\quad P_z\mbox{-a.s.}
\]
for q.e. $z\in E_{0,T}$. Write  $\eta_1=w\wedge v,\, \eta_2=w\vee v$ and denote by $(u,\nu)$ a  solution to OP$(\varphi,f+d\mu,\eta_1,\eta_2)$.
Then $(u(\bX),M,A^\nu)$ is a solution to RBSDE$(\varphi(\bX_{\zeta_\upsilon}), f(\bX,\cdot)+dA^\mu,\eta_1(\bX),\eta_2(\bX))$. Since $(Y^z,M^z,R^z)$ is a solution to the same equation (since  $\eta_1(\bX_t)\le Y^z_t\le \eta_2(\bX_t),\, t\in (0,\zeta_\upsilon)$), we have $(Y^z,M^z,R^z)=(u(\bX),M,A^\nu)$ for q.e. $z\in E_{0,T}$, which implies
that $(u,\nu)$ is a solution to  \mbox{\rm OP}$(\varphi,f+d\mu,h_1,h_2)$.
\end{dow}

\begin{wn}
Under the assumptions of Proposition \ref{stw.red} there exists a unique solution $(u,\nu)$ to \mbox{\rm OP}$(\varphi,f+d\mu,h_1,h_2)$.
\end{wn}

\section{Switching problem}
\label{sec7}

We first describe informally the switching problem. Precise
definitions will be given later on.  Consider a factory in which
we can change a mode of production. Let $c_{j,i}(X)$ be the cost
of the change from mode $j\in\{1,\dots,N\}$ to mode $i$ from some
set $A_j\subset\{1,\dots,j-1,j+1,\dots,N\}$, and let
$f^i(X)+dA^{\mu^i}$ be the payoff rate in mode $i$. Then a
management strategy $\mathcal{S}=(\{\tau_n\},\, \{\xi_n\})$
consists of a pair of two sequences of random variables. The
variable $\tau_n$ is the moment when we decide to switch the mode
of production, and $\xi_n$ is the mode to which we switch at time
$\tau_n$. If $\xi_0=j$ then we start the production at  mode $j$.
If $\xi_0=j$, then under strategy $\mathcal{S}$ the expected
profit on the interval $[0,\zeta_\upsilon]$ is given by the
formula
\begin{align}
\label{eq1.2}
J(z,\mathcal{S},j)&=E_z\Big(\int_0^{\zeta_{\upsilon}}
f^{w^j_r}(\bX_r)\,dr
+\int_0^{\zeta_{\upsilon}}\,dA^{\mu^{w^j_r}}_r\nonumber\\
&\qquad\qquad-\sum_{n\ge 1}
c_{w^j_{\tau_{n-1}},w^j_{\tau_{n}}}(\bX_{\tau_n})
\mathbf{1}_{\{\tau_n<\zeta_\upsilon\}}
+\varphi^{w^j_{\zeta_\upsilon}}(\bX_{\zeta_\upsilon})\Big),
\end{align}
where
\[
w^j_t=j\mathbf{1}_{[0,\tau_1)}(t)+\sum_{n\ge 1} \xi_n
\mathbf{1}_{[\tau_{n},\tau_{n+1})}(t).
\]
The main problem is to find an optimal strategy, i.e. the strategy
$\mathcal{S}^*$ such that
\[
J(x,\mathcal{S}^*,j)=\sup_{\mathcal{S}}J(x,\mathcal{S},j).
\]
In this section we show that $\mathcal{S}^*$ exists and
$\mathcal{S}^*$, $J(z,\mathcal{S}^*,j)$ are determined by a
solution of the system (\ref{eq1.13.1})--(\ref{eq1.13.3}).

\subsection{Systems of BSDEs with oblique reflection}

In what follows $N\in\BN$, $\xi=(\xi^1,\dots,\xi^N)$ is  an
$\FF_T$-measurable random vector, $V=(V^1,\dots,V^N)$ is an
$\mathbb{F}$-adapted process such that $V_0=0$  and each component
$V^j$ is a c\`adl\`ag process of finite variation,
$f:\Omega\times[0,T]\times\BR^N\rightarrow\BR^N$ is a measurable
function such that for every $y\in\BR^N$ the process $f(\cdot,y)$
is $\mathbb{F}$-progressively measurable. Consider a family
$\{h_{j,i}; i,j=1,\dots,N\}$ of measurable functions
$h_{j,i}:\Omega\times[0,T]\times\BR\rightarrow \BR$ such that
$h_{j,i}(\cdot,y^i)$ is progressively measurable for every
$y\in\BR$. For given sets $A_j\subset
\{1,\dots,j-1,j+1,\dots,N\}$, $j=1,\dots,N$, we set
\[
H^j(t,y)=\max_{i\in A_j}h_{j,i}(t,y_i),\quad H(t,y)=(H^1(t,y),
\dots,H^N(t,y)),\quad t\in [0,T],\, y\in\BR^N.
\]

We consider the following system of BSDEs  with oblique reflection:

\begin{equation}
\label{eq8.1} \left\{
\begin{array}{l}
Y^{j}_{t}=\xi^{j}+\int_{t}^{T}f^{j}(r,Y_r)\,dr +\int_t^T
dV_r+\int_t^T dK^{j}_r-\int_{t}^{T}\,dM_r^j,\quad t\in [0,T] ,
\medskip\\
Y^{j}_{t}\ge H^j(t,Y_t),\quad t\in [0,T],
\medskip\\
\int_0^T (Y^{j}_{t-}-[H^j(\cdot,Y)]_{t-})\,dK^{j}_{t}=0, \,j=1,\dots,N.
\end{array}
\right.
\end{equation}

\begin{df}
We say that a triple $(Y,M,K)$ of adapted c\`adl\`ag processes is
a solution of BSDE with oblique reflection (\ref{eq8.1}) if  $Y$
is of class (D), $M$ is a local martingale with $M_0=0$, $K$ is an
increasing process with $K_0=0$ and (\ref{eq8.1}) is satisfied.
\end{df}

If $A_j=\emptyset$, then by convention, $H^j=-\infty$, so $Y^j$
has no lower barrier. We then take $K^j=0$ in the above
definition.

\subsection{Systems of quasi-variational inequalities}

Fix $N\ge1$. Let $\mu^j$,$j=1,\dots,N$, be smooth measures on
$E_{0,T}$, and let $\varphi^j:E\rightarrow\BR$, $f^j:\bar
E_{0,T}\times\BR^N\rightarrow\BR$, $h_{j,i}:\bar E_{0,T}\times\BR
\rightarrow\BR$, $i,j=1,\dots,N$, be measurable functions. We set
\[
f^j(z,y;a)=f^j(z,y_1,\dots,y_{j-1},a,y_{j+1},\dots,y_N),\quad
y\in\BR^N,\, a\in\BR,
\]
and for given sets $A_j\subset \{1,\dots,j-1,j+1,\dots,N\}$,
$j=1,\dots,N$, we set
\[
H^j(z,y)=\max_{i\in A_j}h_{j,i}(z,y_i),\quad H(z,y)=(H^1(z,y),
\dots,H^N(z,y)),\quad z\in\bar E_{0,T},\, y\in\BR^N.
\]
We adopt the convention that the maximum over the empty set equals
$-\infty$. Consequently, if $A_j=\emptyset$ for some $j$, then
$H^j(z,y)=-\infty$.

Set
\[
\varphi=(\varphi^1,\dots,\varphi^N),\quad f=(f^1,\dots,f^N),\quad
\mu=(\mu^1,\dots,\mu^N),
\]
and consider the following system of equations:
\begin{equation}
\label{eq7.3} -\frac{\partial u^j}{\partial
t}-L_tu^j=f^j(t,x,u)+\mu^j\quad\mbox{in }E_{0,T}, \quad
u^j(T,\cdot)=\varphi^j\quad\mbox{on }E
\end{equation}
for $j=1.\dots,N$. In the sequel we denote (\ref{eq7.3}) by
PDE$(\varphi,f+d\mu)$.

\begin{df}
We say that measurable function
$u=(u^1,\dots,u^N):\bar{E}_{0,T}\rightarrow\BR^N$ is a subsolution
(resp. supersolution) of PDE$(\varphi,f+d\mu)$ if there exist
positive measures $\beta^j\in\mathbb{M}$ and $\underline{\varphi}\le \varphi,\, \underline{\varphi}\cdot m\otimes\delta_{\{T\}}\in \mathbb{M}_T$ (resp. $\overline{\varphi}\ge \varphi,\, \overline{\varphi}\cdot m\otimes\delta_{\{T\}}\in \mathbb{M}_T$) such that $u^j$ is a
solution of PDE$(\underline{\varphi}^j,f^j+d\mu^j-d\beta^j)$ (resp.
PDE$(\overline{\varphi}^j,f^j+d\mu^j+d\beta^j$)), $j=1,\dots,N$.
\end{df}

\begin{df}
We say that a quasi-c\`adl\`ag function $u=(u^1,\dots,u^N):\bar
E_{0,T}\rightarrow\BR^N$ is a solution of
(\ref{eq1.13.1})--(\ref{eq1.13.3}) if there exist  positive
smooth measures $\nu^j$ on $E_{0,T}$ such that $(u^j,\nu^j)$ is a solution of
${\rm\underline{O}P}(\varphi^j,f^j(\cdot,u;\cdot)+d\mu^j,H^j(\cdot,u)),\, j=1,\dots, N$.
\end{df}

Let us consider the following hypotheses:
\begin{enumerate}[({A}1)]
\item $\mu, \varphi\cdot m\otimes \delta_{\{T\}}\in\mathbb{M}_T$,
\item for $j=1,\dots,N$ the function is $a\mapsto f^j(z,y;a)$
is nonincreasing
for all $z\in E_{0,T}$, $y\in\BR^N$,
\item $f$ is off-diagonal nondecreasing, i.e. for $j=1,\dots,N$
we have
$f^j(z,y)\le f^j(z,\bar{y})$ for all $y,\bar y\in\BR^N$ such that
$y\le\bar{y}$ and $y^j=\bar{y}^j$,
\item  $y\mapsto f(z,y)$ is continuous for every $z\in E_{0,T}$,
\item  $f^j(\cdot,y)\in qL^1(E_{0,T};m)$ for all $y\in\BR^N$,
$j=1,\dots,N$,
\item there exists a subsolution $\underline{u}$ and a
supersolution $\overline{u}$ of PDE$(\varphi,f+d\mu)$ such that
\[
\underline{u}\le \overline{u},\quad H(\cdot,\overline{u})
\le \overline{u},
\quad \sum_{j=1}^{N}(|f^j(\cdot,\overline{u})|
+|f^j(\cdot,\underline{u};\overline{u}^j)|)\cdot m_1\in\mathbb{M},
\]
\item $H^j$ is continuous  on $\bar{E}_{0,T}\times\BR^N$  with the
product topology consisting of quasi-topology on $\bar{E}_{0,T}$
and Euclidean topology on $\BR^N$, and $h_{j,i}$, $i,j=1,\dots,N$,
are nondecreasing with respect to the second variable.
\end{enumerate}

\begin{tw}
\label{tw8.1} Let the assumptions {{\rm (A1)--(A7)}} hold. Then there
exists a minimal solution  of
\mbox{\rm(\ref{eq1.13.1})--(\ref{eq1.13.3})} such that
$\underline{u}\le u\le \overline{u}$.
\end{tw}
\begin{dow}
First observe that the data $f(\bX,\cdot),H^j(\bX,\cdot),\xi:=
\varphi(\bX_{\zeta_\upsilon}), \underline{Y}:=\underline{u}(\bX),
\overline{Y}:=\overline{u}(\bX)$, $V=A^\mu$ satisfy the assumptions of
\cite[Theorem 3.11]{K:arx} under the measure $P_z$ for q.e. $z\in
E_{0,T}$. Set $u_0=\underline{u}$ and $Y^0=\underline{Y}$. By
Theorem \ref{tw6.2} (see also Remark \ref{uw.main}), for every $n\ge 1$,
\[
u^j_n(X_t)=Y^{n,j}_t,\quad A^{\nu_n}_t=K^{n,j}_t,
\]
where $(u^j_n,\nu^j_n)$ is a solution of
${\rm\underline{O}P}(\varphi^j,f^j(\cdot,u_{n-1};\cdot)+\mu^j, H^j(\cdot,u_{n-1}))$ and
the triple $(Y^{n,j},K^{n,j},M^{n,j})$ is a solution of {\rm
\underline{R}BSDE}$(\xi^j,
f^j(\bX,Y^{n-1};\cdot)+dV^j,H^j(\cdot,Y^{n-1}))$. By
Proposition \ref{stw.cr}, $u_n\le u_{n+1}$ q.e. Set
$u:=\sup_{n\ge 0} u_n$. By \cite[Theorem
3.11]{K:arx},
\[
Y^{n,j}_t\nearrow Y_t,\quad t\in [0,\zeta_\upsilon],\quad P_z\mbox{-a.s.}
\]
for q.e. $z\in E_{0,T}$, where $(Y,M,K)$ is the minimal solution of (\ref{eq8.1}) such that
$\underline{Y}\le Y\le \overline{Y}$.
Hence, since the exceptional sets coincide with the sets of zero capacity, we have
\[
u^j(\bX_t)=Y^j_t,\quad t\in [0,\zeta_\upsilon],\quad P_z\mbox{a.s.}
\]
for q.e. $z\in E_{0,T}$.
We see that the triple $(Y^j,M^j,K^j)$  is a solution to the problem
{\rm\underline{R}BSDE}$(\varphi^j(\bX_{\zeta_\upsilon}),f^j(\bX,u(\bX);\cdot)
+dA^{\mu^j},H^j(\cdot,u(\bX)))$. By Theorem \ref{tw6.2}, $K^j=A^{\nu^j}$, where $(u^j,\nu^j)$ is a solution to
${\rm\underline{O}P}(\varphi^j,f^j(\cdot,u;\cdot)+d\mu^j,H^j(\cdot,u)),\, j=1,\dots,N$,
which implies that the pair $(u,\nu)$ is a solution of (\ref{eq1.13.1})--(\ref{eq1.13.3}).
Minimality of $u$ follows from the  minimality of $Y$.
\end{dow}

\begin{uw}
\label{uw4.4}
Let $u$ be the minimal solution of Theorem \ref{tw8.1}. Observe that under assumptions of Theorem \ref{tw8.1} $f^j(\cdot,u)\cdot m_1\in \mathbb M$
and $\nu^j\in\mathbb M,\, j=1,\dots,N$. Indeed, first observe that $\bar u^j\ge H^j(\cdot,u)$ and
\begin{equation}
\label{eq.n121}
f^j(\cdot,\bar u)\le f^j(\cdot,u;\bar u)\le f^j(\cdot,\underline u;\bar u).
\end{equation}
Since $u^j$ is a solution of ${\rm\underline{O}P}(\varphi^j,f^j(\cdot,u;\cdot)+d\mu^j,H^j(\cdot,u)),\, j=1,\dots,N$
we get the result by Remark \ref{uw.main}.
Moreover, if we assume that  $\EE_\gamma$ has the dual Markov
property for some  $\gamma\ge 0$,  assumptions (A1), (A6),
are satisfied with $\BM$ replaced by $\MM_1$ and
the measures $\beta^j$ appearing in the definition of the
supersolution $ \overline{u}$ belong to $\MM_1$.  Then by Remark \ref{uw.main} under assumptions of Theorem \ref{tw8.1}
$\nu^j\in\MM_1$, and $f^j(\cdot,u)\in L^1(E_{0,T};m_1)$, $j=1,\dots,N$.
\end{uw}

Let us consider the following hypothesis:
\begin{enumerate}
\item [(A8)] there exists a subsolution $\underline{u}$ and a
supersolution $\overline{u}$ of PDE$(\varphi,f+d\mu)$ and  a function $v=(v^1,\dots,v^N)$ which is a difference of potentials on $\bar E_{0,T}$ such that
\[
\sum_{j=1}^N(|f^j(\cdot,\overline{u}; v^j)|
+|f^j(\cdot,\underline{u}; v^j)|)\cdot m_1\in \mathbb{M}
.\]
\end{enumerate}

\begin{stw}
\label{stw8.1} Let assumptions  \mbox{\rm (A1)--(A5), (A8)} hold.
Then there exists  minimal  solution $u$ of \mbox{\rm
PDE}$(\varphi,f+d\mu)$ such that $\underline{u}\le u\le
\overline{u}$ q.e.
\end{stw}
\begin{dow}
Observe that $f(\bX,\cdot)$, $\underline{Y}:= \underline{u}(\bX)$,
$\overline{Y}:=\overline{u}(\bX)$, $S=v(\bX)$ (see Remark \ref{uw6.1}), $V=A^\mu$,
$\xi:=\varphi(\bX_{\zeta_\upsilon})$ satisfy the assumptions of
\cite[Theorem 2.12]{K:arx} under the measure $P_z$ for q.e. $z\in
E_{0,T}$. Set $u_0=\underline{u}$. By \cite[Theorem 5.8]{K:JFA},
$Y^{n,j}=u^j_n(\bX)$, where $u^j_n$ is a solution of
PDE$(\varphi,f^j(\cdot,u_{n-1};\cdot)+d\mu^j)$ and $(Y^{n,j},
M^{n,j})$ is a solution of
BSDE$(\xi^j,f^j(\bX,Y^{n-1};\cdot)+dV^j),\, j=1,\dots,N$.  By \cite[Corollary 5.9]{K:JFA}  $u_n\le u_{n+1}$ q.e. Set $u=\sup_{n\ge 1} u_n$.
By \cite[Theorem 2.12]{K:arx} it follows that
$Y^n_t\nearrow Y_t,\, t\in [0,\zeta_\upsilon],\, P_z$-a.s. for q.e. $z\in E_{0,T}$, where $(Y,M)$ is a  minimal solution of
BSDE$(\varphi(\bX_{\zeta_\upsilon}), f(\bX,\cdot)+dA^\mu)$ such
that $\underline{u}(\bX)\le Y\le\overline{u}(\bX)$.  Hence (since the exceptional sets coincide  with the sets of  zero capacity)  $u(\bX_t)=Y_t,\, t\in [0,\zeta_\upsilon],\, P_z$-a.s. for q.e. $z\in E_{0,T}$, which implies that $u$ is
the minimal solution of PDE$(\varphi, f+d\mu)$ such that
$\underline{u}\le u\le \overline{u}$ q.e.
\end{dow}

\begin{tw}
Assume {{\rm (A1)--(A7)}}. Then there exists  minimal solution
$u_n$  of the system
\begin{equation}
\label{eq8.2}
-\frac{\partial u_n}{\partial t}-L_t u_n^j
=f^j(\cdot,u_n)+n(u_n^j-H^j(\cdot,u_n))^-+\mu
\end{equation}
such that $\underline{u}\le u_n\le\overline{u}$. Moreover,
$u_n\nearrow u$ q.e., where $u$ is  minimal solution of \mbox{\rm
(\ref{eq1.13.1})--(\ref{eq1.13.3})} such that $\underline{u}\le
u\le \overline{u}$.
\end{tw}
\begin{dow}
Observe that $\overline{u}$ is a supersolution of (\ref{eq8.2}),
and $\underline{u}$ is a subsolution of (\ref{eq8.2}). Moreover,
(A8) for (\ref{eq8.2}) is satisfied with $v=\overline{u}$. By Proposition
\ref{stw8.1} there exists  a minimal solution $u_n$ of
(\ref{eq8.2}). By the definition and construction of minimal solution to (\ref{eq8.2}) (see Proposition \ref{stw8.1}) and minimal solution to BSDE$(\varphi(\bX_{\zeta_\upsilon}),f_n(\bX,\cdot)+dA^\mu)$  (see
\cite[Theorem 2.12]{K:arx}), $u_n(\bX)$ is the
first component of minimal solution of
BSDE$(\varphi(\bX_{\zeta_\upsilon}),f_n(\bX,\cdot)+dA^\mu)$ with
$f^j_n(z,y)=f^j(z,y)+n(y^j-H^j(z,y))^-$. By \cite[Theorem
3.15]{K:arx}, the  sequence $\{u_n(\bX)\}$ is nondecreasing and
$u_n(\bX)_t\nearrow Y_t,\, t\in [0,\zeta_\upsilon],\, P_z$-a.s. for q.e. $z\in
E_{0,T}$, where $Y$ is the first component of the minimal solution $(Y,M,K)$ of
(\ref{eq8.1}) such that $\underline{u}(\bX)\le
Y\le \overline{u}(\bX)$. Since the  sets coincide  with the sets of  zero capacity $u_n\le u_{n+1}$, q.e. on $E_{0,T}$. Let $u:=\sup_{n\ge 1}u_n$. It is clear
that $Y_t=u(\bX_t),\, t\in [0,\zeta_\upsilon],\, P_z$-a.s. for q.e. $z\in E_{0,T}$.
Now we see that $(Y^j,M^j,K^j)$ is a solution to {\rm
\underline{R}BSDE}$(\varphi^j(\bX_{\zeta_\upsilon}),
f^j(\bX,u(\bX);\cdot)+dA^{\mu^j},H^j(\cdot,u(\bX)))$. By Theorem \ref{tw6.2} $K^j=A^{\nu^j}$, where
$(u^j,\nu^j)$ is a solution to ${\rm\underline{O}P}(\varphi^j,f^j(\cdot,u;\cdot)+d\mu^j,H^j(\cdot,u)),\, j=1,\dots,N$.
This implies that $(u,\nu)$ is the minimal solution of
(\ref{eq1.13.1})--(\ref{eq1.13.3}) such that $\underline{u}\le
u\le \overline{u}$. Of course,  $u_n\nearrow u$
q.e.
\end{dow}

\subsection{Value function for the switching problem }

In what follows we assume that $H^j$ are of the form
\begin{equation}
\label{eq8.3}
H^j(z,y)=\max_{i\in A_j}(-c_{j,i}(z)+y^i),
\end{equation}
where $c_{j,i}$ are  quasi-continuous functions on $E_{0,T}$ such
that for some  constant $c>0$,
\[
c_{j,i}(z)\ge c,\quad z\in E_{0,T},\quad i\in A_j,\quad
j=1,\dots,N.
\]

By a strategy we call a pair $\mathcal{S}=(\{\xi_n\}, \{\tau_n\})$
consisting of a sequence $\{\tau_n,n\ge 1\}$ of increasing
$\BF$-stopping times such that
\[
P_z(\tau_n<\zeta_\upsilon,\,\forall\,\, n\ge 1)=0
\]
for q.e. $z\in E_{0,T}$, and a sequence $\{\xi_n,n\ge1\}$ of
random variables taking values in $\{1,\dots, N\}$ such that
$\xi_n$ is $\FF_{\tau_n}$-measurable for each $n\ge1$. The set of
all strategies we denote by $\mathbf{A}$. For
$\mathcal{S}\in\mathbf{A}$ we set
\[
w^j_t=j\mathbf{1}_{[0,\tau_1)}(t)+\sum_{n\ge 1}\xi_n
\mathbf{1}_{[\tau_{n},\tau_{n+1})}(t).
\]

\begin{uw}
\label{uw3.3} In Theorem \ref{stw8.1} assume additionally that
$\mu\in\mathbb{M}_c$, $h_{j,i}$ are strictly increasing  with
respect to $y$, and that the following condition considered in
\cite{HZ} is satisfied:
\begin{enumerate}
\item[(A9)] there are no $(y_1,\dots, y_k)\in\BR^k$ and
$j_2\in A_{j_1},\dots,j_k\in A_{j_{k-1}}, j_1\in A_{j_k}$ such
that
\[
y_1=h_{j_1,j_2}(z,y_2), y_2=h_{j_2,j_3}(z,y_3),\dots,
y_{k-1}=h_{j_{k-1},j_k}(z,y_k), y_k=h_{j_k,j_1}(z,y_1).
\]
\end{enumerate}
Then $\nu\in\mathbb{M}_c$ and $u$ is quasi-continuous. This
follows from \cite[Remark 3.14]{K:arx}. Observe that (A9) is
satisfied  for $h_{j,i}$  defined by (\ref{eq8.3}).
\end{uw}

\begin{tw}
Assume that $f$ does not depend on $y$, the functions  $H^j$ are
of the form \mbox{\rm(\ref{eq8.3})} and  $f^j\cdot m_1,\mu^j\in\BM$,
$j=1,\dots,N$. Then there exists a unique solution $u$ of
\mbox{\rm(\ref{eq1.13.1})--(\ref{eq1.13.3})}. Moreover,
\begin{align*}
u^j(z)=\sup_{\mathcal{S}\in\mathbf{A}}J(z,\mathcal{S},j)
\end{align*}
and
\begin{align*}
u^j(z)&=E_z\Big(\int_0^{\zeta_\upsilon}
f^{w^{j,*}_r}(\bX_r)\,dr+\int_0^{\zeta_\upsilon}\,dA^{\mu^{w^{j,*}_r}}_r\\
&\qquad\quad-\sum_{n\ge1} c_{w^{j,*}_{\tau_{n-1}},
w^{j,*}_{\tau_{n}}}(\bX_{\tau_n})
\mathbf{1}_{\{\tau_n<{\zeta_\upsilon}\}}
+\varphi^{w^{j,*}_{\zeta_\upsilon}}(\mathbf{X}_{\zeta_\upsilon})\Big),
\end{align*}
where
\[
w^{j,*}_t=j\mathbf{1}_{[0,\tau^{j,*}_1)}(t)+\sum_{n\ge 1}\xi^{j,*}_n
\mathbf{1}_{[\tau^{j,*}_{n},\tau^{j,*}_{n+1})}(t)
\]
and
\[
\tau^{j,*}_{0}=0,\quad \quad \xi^{j,*}_{0}=j,
\]
\[
\tau_k^{j,*}=\inf\{t\ge \tau^{j,*}_{k-1}:
u^{\xi^{j,*}_{k-1}}(\bX_t)=H^{\xi^{j,*}_{k-1}}(\bX_t,u(\bX_t))\}
\wedge \zeta_\upsilon,\quad k\ge1,
\]
\[
\xi^{j,*}_k=\max\{
i\in A_{\xi^{j,*}_{k-1}};\, H^{\xi^{j,*}_{k-1}}(\bX_{\tau_k^{j,*}},u(\bX_{\tau_k^{j,*}}))
=-c_{\xi^{j,*}_{k-1},i}(\bX_{\tau_k^{j,*}})+u^i({\bX_{\tau_k^{j,*}})}\},\quad
k\ge 1.
\]

\end{tw}
\begin{dow}
Follows from Proposition \ref{stw8.1}, Remark \ref{uw3.3} and
\cite[Theorem 4.3]{K:arx}.
\end{dow}

\end{document}